%% file: main.tex
\documentclass[11pt]{article}
\usepackage[margin=1in]{geometry}
\usepackage[T1]{fontenc}
\usepackage{times}
\usepackage[round]{natbib}
\usepackage{parskip}

\input{math_commands.tex}

\usepackage[colorlinks=true,linkcolor=blue!50!black,citecolor=blue!50!black,urlcolor=blue!50!black]{hyperref}
\usepackage{url}
\usepackage{algorithm}
\usepackage{algpseudocode}
\usepackage{tikz}
\usetikzlibrary{shapes.geometric, arrows.meta, positioning, fit, backgrounds}
\usepackage{subcaption}
\usepackage{graphicx}
\usepackage{booktabs}
\usepackage{multirow}
\usepackage{makecell}
\usepackage{longtable}
\usepackage{authblk}
\usepackage{placeins}
\newcommand{\rev}[1]{#1}

\newtheorem{proposition}{Proposition}
\newtheorem{proof}{Proof}

\title{Beyond Decoupled Decisions: Jointly Learning Label Extension and Column Selection in Column Generation for Vehicle Routing Problems}
\author[1]{Hongli Li}
\author[1]{Zengxiang Lei}
\author[2]{Xinwu Qian}
\author[1,*]{Satish V. Ukkusuri}
\affil[1]{Lyles School of Civil and Construction Engineering, Purdue University, West Lafayette, IN, USA}
\affil[2]{Department of Civil and Environmental Engineering, Rice University, Houston, TX 77005, USA}
\affil[*]{Corresponding author}
\affil[ ]{\small ORCID: Hongli Li (\href{https://orcid.org/0009-0001-7443-5563}{0009-0001-7443-5563}), Zengxiang Lei (\href{https://orcid.org/0000-0002-7639-438X}{0000-0002-7639-438X}), Xinwu Qian (\href{https://orcid.org/0000-0001-8001-2164}{0000-0001-8001-2164}), Satish V. Ukkusuri (\href{https://orcid.org/0000-0001-8754-9925}{0000-0001-8754-9925})}
\date{}

\begin{document}
\maketitle

\begin{abstract}
Column generation (CG) is a powerful framework for large-scale routing optimization that avoids exhaustive enumeration by iteratively solving master and pricing problems. However, its efficiency is constrained by two coupled bottlenecks: the labeling-based pricing problem spends substantial computation exploring labels and uncovering many redundant candidate columns, while only a small, complementary subset materially improves the master problem. We study CG with routing pricing subproblems solved by labeling algorithms and formulate label-extension ordering and column selection as two sequential decisions within a shared Markov decision process. We consider a practical setting that learns from rich pre-collected CG trajectories without solving CG online during training, using data generated by existing CG solvers. We propose \textsc{Coupled-CG}, a two-stage offline Q-learning framework with stage-specific conservative regularization that combines support-aware value learning for extension ordering, relation-conditioned marginal column selection, and coupled Bellman updates that propagate value across the two stages and successive CG iterations. Experiments on public CVRPTW benchmarks demonstrate substantial reductions in labeling effort and end-to-end CG runtime. Ablation studies quantify the additional gains from jointly learning the two decisions, while instance-level analyses identify the structural characteristics of instances that benefit most from Coupled CG.
\end{abstract}

\input{Content/Intro}

\input{Content/Literature}

\input{Content/Problem}

\input{Content/Methodology}
\input{Content/Experiment}

\input{Content/Conclusion}

\subsection*{Funding}
This work was supported by the National Center for Transportation Cybersecurity and Resiliency (TraCR), a U.S. Department of Transportation National University Transportation Center headquartered at Clemson University, and by the U.S. Department of Energy Co-E3T project (Energy-Efficient and Equitable Transit through User-Centric Hardware and Software Co-Development and Community Co-Design).

\subsection*{AI use statement}
In this work, we used ChatGPT to assist with writing and language refinement, including improving clarity, grammar, and presentation of the manuscript. We also used ChatGPT during implementation to generate preliminary code drafts and suggestions for selected components. All AI-assisted text was reviewed and revised by the authors, and all AI-generated or AI-assisted code was inspected, tested, and modified as needed to ensure correctness and consistency with the proposed methodology. ChatGPT was not used to generate experimental results, perform evaluations, or make final scientific conclusions. We take responsibility for the final content of this work, including all text, claims, code, and other artifacts produced with the aid of generative AI.

\subsection*{Ethics statement}
This work does not involve human subjects, personally identifiable information, or sensitive private data. The study uses public benchmark instances and focuses on algorithmic methods for accelerating column generation. We are not aware of any specific ethical concerns arising from the methodology or experiments presented in this work.

\bibliographystyle{plainnat}
\bibliography{references}

\input{Content/Appendix}
\end{document}

%% file: math_commands.tex
\usepackage{amsmath,amsfonts,bm}

\def\eqref#1{equation~\ref{#1}}

\def\1{\bm{1}}

\DeclareMathAlphabet{\mathsfit}{\encodingdefault}{\sfdefault}{m}{sl}
\SetMathAlphabet{\mathsfit}{bold}{\encodingdefault}{\sfdefault}{bx}{n}



%% file: Content/Intro.tex
\section{Introduction}

Column generation (CG) solves large-scale linear programs by generating decision variables (columns) on demand~\citep{desaulniers2006column}. At each iteration, a restricted master problem (RMP) optimizes over the current columns, while a pricing subproblem searches for new columns that can improve the solution. 
In routing problems, a column typically represents a feasible route, and the pricing problem becomes seeking negative-reduced-cost routes under operational constraints such as vehicle capacity and customer time windows~\citep{irnich2005shortest, Costa2019}.
This decomposition has substantially advanced exact routing optimization, enabling, for example, capacitated vehicle routing instances with up to 360 customers that had previously been tackled only heuristically~\citep{pecin2017improved}. We focus on this important class of CG methods in which the pricing subproblem is a constrained routing problem~\citep{righini2008new,Baldacci2011,dabia2013branch,he2019branch,li2022two,li2026anticipatory,mhamedi2026branch}.

Within this class, we consider pricing methods based on labeling algorithms, which construct routes incrementally while tracking resources such as capacity and time and pruning infeasible or dominated partial routes~\citep{Feillet2004,irnich2005shortest}. Despite their effectiveness, these CG methods face two interacting computational bottlenecks. First, pricing remains a generally NP-hard constrained routing problem. Because the number of partial routes (labels) can grow exponentially, the labeling algorithm needs to explore and compare a large number of labels at every CG iteration. Second, after candidate columns are generated, selecting which ones to add to the RMP remains challenging. Reduced cost provides a natural measure of immediate improvement, but a column with a more negative reduced cost does not necessarily yield faster overall convergence to CG termination~\citep{chi2022deep}.
Adding many columns can reduce the number of CG iterations, yet it also enlarges the RMP and increases its solution cost~\citep{Morabit2021}. 

These two bottlenecks are coupled across CG iterations.  Label-extension decisions determine both pricing effort and the columns available to the RMP, while the column-selection decisions alter the RMP solution and the dual information that guides subsequent pricing. \rev{This coupling is part of the CG algorithm itself, not a modeling assumption.} 
Together, these challenges limit the scale of routing problems that can be solved exactly~\citep{desaulniers2019selective,cabrera2023solving}.

Existing CG acceleration studies address pricing and column selection separately. On the pricing side, classical methods reduce labeling effort through stronger dominance and bounding rules, bidirectional search, and route relaxations~\citep{Righini2006,Tilk2017,desaulniers2019selective,Sadykov2021}; more recent learning-based approaches further reduce the pricing search space through learned arc or heuristic selection~\citep{Xu2023RLHH,morabit2023,Xia2024NeuralCG}. For column selection, classical and learning-based strategies selectively admit generated columns to control RMP growth and accelerate CG termination 
~\citep{touati2010solutions,Morabit2021,chi2022deep,Yuan2024}. These approaches can yield substantial computational gains, including reported runtime reductions of up to 40\% from learned pricing guidance~\citep{morabit2023} and 30\% from learned column selection~\citep{Morabit2021}. \rev{However, each of these methods learns one decision and fixes the other. We hypothesize that a state-dependent joint policy for extension ordering and column selection requires less cumulative labeling effort than fixed rules or independently learned policies, and that such a policy can be learned from logged CG runs.}

In this study, we address this gap by formulating CG as a Markov decision process (MDP) with two sequentially coupled policies: one selects a label-extension ordering rule for pricing, while the other selects a subset of the resulting columns for the RMP. \rev{Unlike learned pricing methods that prune or replace the pricing network, the pricing policy changes only the order in which labels are extended. The feasibility and dominance rules and the final exact pricing step are unchanged, so CG still terminates at the LP optimum.} We design stage-specific rewards that penalize labeling effort and measure RMP objective improvement; the latter telescopes over an episode, providing intermediate feedback without altering our surrogate CG acceleration objective of minimizing cumulative labeling effort. 
\rev{An online training episode would be a complete CG run, with labeling and an LP re-solve at every iteration. We therefore learn offline from a fixed set of trajectories logged from CG runs on training instances.} We develop \textsc{Coupled-CG}, a two-stage offline conservative Q-learning framework. It uses support-aware conservative learning to address limited action coverage and relation-aware marginal column selection to capture interactions among candidate columns. The two stages are further linked through coupled Bellman updates: the pricing-stage target incorporates downstream column-selection value, while the column-selection target incorporates the next iteration's pricing value. This enables both policies to account for their longer-term effects across CG iterations.
We evaluate \textsc{Coupled-CG} on public benchmarks of the \rev{capacitated vehicle routing problem with time windows (CVRPTW)}~\citep{solomon1987algorithms}. Our contributions are:

\begin{enumerate}
    \item \rev{\textbf{Formulation}. We cast extension ordering and column selection as two stages of one MDP, with stage rewards that preserve the cumulative labeling-effort objective (Proposition~\ref{prop:reward_consistency}).}
    \item \rev{\textbf{Coupled offline learning}. Cross-stage Bellman updates value an ordering rule by the columns it yields and a column set by the pricing problem it creates. Support-aware ordering and relation-conditioned selection make the policy learnable from logged data; removing them raises mean labels by 73\% and 432\%.}
    \item \rev{\textbf{Evidence}. On Solomon CVRPTW instances, \textsc{Coupled-CG} reduces CG time by 39--77\% and explored labels by 44--78\% relative to Plain CG, improves on independently trained and hybrid variants, and extends to larger instances and two other problems.}
\end{enumerate}

%% file: Content/Literature.tex
\section{Related Work}\label{sec:literature_view}
\paragraph{Pricing and Labeling Acceleration}

Classical methods accelerate routing pricing through route relaxations, bidirectional search, stronger dominance and bounding rules, specialized graph representations, and selective pricing~\citep{Baldacci2011,Righini2006,Tilk2017,Sadykov2021,Enerback2024,desaulniers2019selective}. Learning-based approaches further guide pricing by reducing the search network through arc selection~\citep{morabit2023}, guiding or selecting pricing heuristics~\citep{Shen2022MLPH,Xu2023RLHH}, and identifying or generating promising candidate columns~\citep{Xia2024NeuralCG,Abouelrous2025RLPricing}. These approaches improve route generation, but pricing efficiency alone does not capture its impact on overall CG performance. We instead learn extension ordering jointly with downstream column selection to account for their interaction across CG iterations.

\paragraph{Column Selection and Coordination with Pricing}
Classical column-selection methods use reduced-cost criteria and diversification to control RMP growth~\citep{touati2010solutions}, while learning-based approaches predict column usefulness~\citep{Morabit2021}, model selection sequentially~\citep{chi2022deep}, or select multiple columns per iteration~\citep{Yuan2024}. FFCG~\citep{hu2025ffcg} further performs sequential marginal-value selection without enumerating subsets. Our contribution is not sequential selection itself, but its integration with learned pricing-order control in a coupled offline framework, using relation-conditioned set values and cross-stage credit assignment. Existing pricing and column-selection components have also been combined~\citep{Morabit2021,morabit2023}, but are trained separately; we instead learn the two decisions jointly as sequentially coupled policies.

\paragraph{Offline Reinforcement Learning (RL)}

Offline RL learns policies from fixed datasets without further environment interaction during training~\citep{prudencio2023survey}, while conservative value learning mitigates overestimation under limited action coverage~\citep{fujimoto2019off,kumar2020conservative}. We adapt these ideas to the coupled structure of CG, where pricing and column selection have distinct action spaces but sequentially affect one another. \textsc{Coupled-CG} combines stage-specific conservative learning with relation-aware column selection and cross-stage Bellman targets that propagate value between the two decisions across CG iterations.

%% file: Content/Problem.tex
\section{Problem Formulation} \label{sec:problem_setup}

\subsection{Preliminary}\label{sec:preliminary}
We use the CVRPTW as a representative setting to introduce the CG framework considered in this work. Let $\mathcal N$ be the customer set and $\mathcal R$ the set of feasible routes, where each route starts and ends at the depot, visits each customer at most once, and satisfies capacity and time-window constraints. Route $r\in\mathcal R$ has cost $c_r$, and $\alpha_{i,r}\in\{0,1\}$ indicates whether it serves customer $i$. With fleet limit $V$ and route-selection variable $X_r$, the route-based formulation is
\begin{equation}
    \mathcal M_{\text{IP}} = \min\; \sum_{r\in\mathcal R} c_r X_r
\quad \text{s.t.}\quad
\sum_{r\in\mathcal R}\alpha_{i,r}X_r=1,\ \forall i\in\mathcal N,\quad
\sum_{r\in\mathcal R}X_r\le V,\quad
X_r\in\{0,1\},\ \forall r\in\mathcal R,
\end{equation}

Because $|\mathcal R|$ can grow exponentially with $|\mathcal N|$, CG solves the LP relaxation using only a restricted route set $\mathcal R'\subseteq\mathcal R$. Solving the resulting restricted master problem (RMP) yields dual variables $d_i$ for the customer-coverage constraints and $q$ for the fleet constraint, defining the reduced cost of route $r$ as
$\bar{c}_r=c_r-\sum_{i \in \mathcal N} d_i \alpha_{i,r}-q$. The pricing problem then searches $\mathcal R$ for negative-reduced-cost routes. Selected routes are added to $\mathcal R'$, the RMP is re-solved, and the process repeats until no improving route remains, at which point the full LP relaxation is solved. CG can subsequently be embedded within branch-and-bound for exact integer optimization. We focus on accelerating the CG loop; notation is summarized in Appendix~\ref{appendix:notation_tab}.

\subsection{Problem Formulation: Coupled Acceleration of Column Generation}
\label{sec:theory}

\paragraph{Acceleration objective.}

We seek a CG acceleration strategy that generalizes across problem instances $\mathcal I\sim\mathcal D$, where $\mathcal D$ represents the instance distribution of the target application and may be specified through an instance sampler. For each instance, CG repeatedly solves a pricing problem $P_t$  and reoptimizes the RMP until termination (Section~\ref{appendix:CG}). Acceleration therefore requires considering both the computational effort of individual pricing calls and the progress induced by the generated and selected columns.

\paragraph{Controlling pricing search through extension ordering.}

The pricing problem is first solved by heuristic labeling with early stopping, which terminates when the number of explored labels reaches a predefined limit $\Gamma$; exact labeling is invoked when needed for optimality certification. At iteration $t$, an extension-ordering rule $g_t\in\mathcal G$ determines the order in which active labels are processed in pricing problem $P_t$, \rev{where $\mathcal G$ denotes the set of five candidate rules based on minimum reduced cost, minimum or maximum accumulated time, and minimum or maximum accumulated load.} We consider strategies that modify only this processing order while leaving the extension, feasibility, reduced-cost, and valid dominance rules unchanged. Details of both labeling procedures are provided in Section~\ref{appendix:labeling_algorithm}.

The effect of extension ordering differs between exact and early-stop labeling. Under exact labeling, ordering can change the search trajectory and computational effort, for example, by discovering dominating labels earlier and avoiding unnecessary extensions, but it does not change the pricing optimum: all orderings solve the same pricing problem exactly and therefore return the same optimal reduced-cost value. Under early stopping, however, the ordering also determines which parts of the search space are explored within the fixed budget $\Gamma$, and can therefore change the candidate routes and reduced-cost values returned. Thus, extension ordering affects pricing effort in both settings and, under early stopping, additionally shapes the candidate columns available to the subsequent RMP update.

\paragraph{Joint pricing and column-selection objective.}

At iteration $t$, extension ordering shapes the candidate set $\mathcal C_t$, while the selected subset $O_t\subseteq\mathcal C_t$ determines the RMP update and the dual information defining the next pricing problem:
\begin{equation}
g_t \longrightarrow
\mathcal C_t \longrightarrow
O_t \longrightarrow
\text{RMP update} \longrightarrow P_{t+1}.
\end{equation}
Pricing decisions therefore affect the choices available to column selection, while column-selection decisions affect subsequent pricing effort. We optimize these decisions jointly using cumulative labeling effort as a surrogate for CG computation. Let $N_t(g_t)$ denote the number of labels explored in pricing call $t$, including heuristic labeling and exact fallback when invoked. With extension-ordering policy $\pi^P$ and column-selection policy $\pi^C$, our objective is
\begin{equation}
\mathcal M:
\min_{\pi^P,\pi^C}
\mathbb E_{\mathcal I\sim\mathcal D}
\left[
\sum_{t=1}^{T(\pi^P,\pi^C;\mathcal I)}
N_t(g_t)
\right],
\end{equation}
where $T(\pi^P,\pi^C;\mathcal I)$ denotes the number of CG iterations until exact pricing certifies termination. This objective captures both per-call labeling effort and the effects of the two policies on subsequent pricing problems and trajectory length.

\subsection{Column Generation Framework}\label{appendix:CG}
\begin{figure}[!htbp]
\centering

\resizebox{\columnwidth}{!}{%
\begin{tikzpicture}[
    >=Latex,
    node distance=1.8mm,
    every node/.style={font=\scriptsize},
    block/.style={
        rectangle,
        rounded corners,
        draw,
        align=center,
        minimum width=1.35cm,
        minimum height=6.5mm,
        inner sep=1.5pt
    },
    decision/.style={
        diamond,
        draw,
        align=center,
        aspect=1.55,
        text width=1.25cm,
        inner sep=0.3pt
    },
    rulelearner/.style={
        rectangle,
        rounded corners,
        draw=blue!70!black,
        fill=blue!12,
        thick,
        align=center,
        minimum width=1.35cm,
        minimum height=6.5mm,
        inner sep=1.5pt
    },
    columnlearner/.style={
        rectangle,
        rounded corners,
        draw=orange!80!black,
        fill=orange!15,
        thick,
        align=center,
        minimum width=1.35cm,
        minimum height=6.5mm,
        inner sep=1.5pt
    },
    line/.style={->, thick}
]

\node[block] (rmp)
{Solve RMP};

\node[rulelearner, right=of rmp] (ruleH)
{Select\\extension rule};

\node[block, right=of ruleH] (heur)
{Heuristic\\labeling};

\node[decision, right=of heur] (heurcheck)
{Negative\\columns\\found?};

\node[rulelearner, right=of heurcheck] (ruleE)
{Select\\extension rule};

\node[block, right=of ruleE] (exact)
{Exact\\labeling};

\node[decision, right=of exact] (exactcheck)
{Negative\\columns\\found?};

\node[block, right=of exactcheck] (term)
{Terminate};

\node[
    columnlearner,
    above=7mm of heurcheck
] (selectH)
{Select\\columns};

\node[
    columnlearner,
    below=7mm of exactcheck
] (selectE)
{Select\\columns};

\draw[line] (rmp) -- (ruleH);

\draw[line] (ruleH) -- (heur);

\draw[line] (heur) -- (heurcheck);

\draw[line]
    (heurcheck) --
    node[above, font=\scriptsize] {No}
    (ruleE);

\draw[line] (ruleE) -- (exact);

\draw[line] (exact) -- (exactcheck);

\draw[line]
    (exactcheck) --
    node[above, font=\scriptsize] {No}
    (term);

\draw[line]
    (heurcheck.north) --
    node[right, font=\scriptsize] {Yes}
    (selectH.south);

\draw[line]
    (selectH.west)
    -|
    node[
        pos=0.30,
        above,
        font=\scriptsize
    ] {add selected columns}
    (rmp.north);

\draw[line]
    (exactcheck.south) --
    node[right, font=\scriptsize] {Yes}
    (selectE.north);

\draw[line]
    (selectE.west)
    -|
    node[
        pos=0.30,
        below,
        font=\scriptsize
    ] {add selected columns}
    (rmp.south);

\end{tikzpicture}%
}

\caption{
\textsc{Coupled-CG} framework. The extension-ordering learner guides heuristic
and exact labeling, while the column-selection learner selects columns
to update the RMP.
}
\label{fig:coupled_cg}

\end{figure}
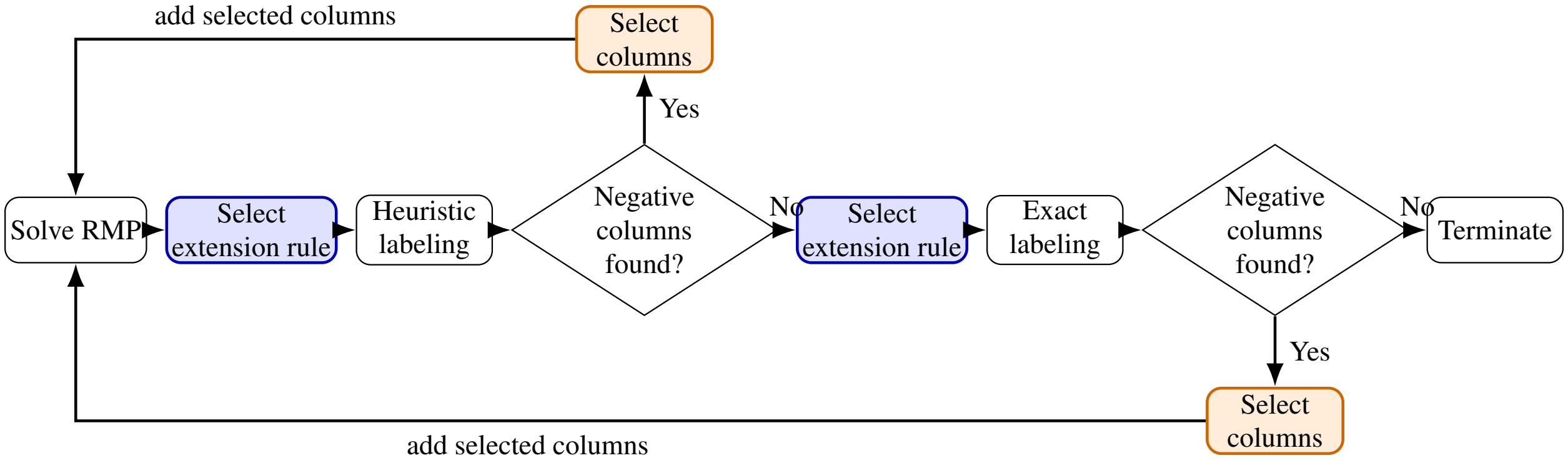

Figures~\ref{fig:original_cg} and~\ref{fig:coupled_cg}
illustrate the standard CG procedure and the proposed
\textsc{Coupled-CG} framework, respectively. In the original CG
framework, each iteration first solves the RMP and then
invokes heuristic labeling to quickly search for negative
reduced-cost columns. If heuristic labeling fails to identify
any valid improving column, exact labeling is used as a fallback. Unlike heuristic labeling, which may terminate early and therefore does not guarantee pricing optimality, exact labeling solves the pricing problem exactly and certifies whether any negative reduced-cost column remains. CG terminates only when exact labeling finds no improving column.

\textsc{Coupled-CG} preserves this standard pricing structure while
introducing learning at two decision points. Before both heuristic and exact labeling, an extension-ordering learner selects the rule used to guide label processing. After valid columns are generated, a column-selection learner selects the subset to be added to the RMP. Thus, the learned components control both pricing search and the RMP update, while exact labeling remains the fallback mechanism that preserves the standard CG termination condition.

\begin{figure}[!htbp]
\centering
\begin{tikzpicture}[
    >=Latex,
    node distance=3.5mm,
    every node/.style={font=\scriptsize},
    block/.style={
        rectangle,
        rounded corners,
        draw,
        align=center,
        minimum width=1.75cm,
        minimum height=7mm,
        inner sep=2pt
    },
    decision/.style={
        diamond,
        draw,
        align=center,
        aspect=1.6,
        text width=1.55cm,
        inner sep=0.5pt
    },
    line/.style={->, thick}
]

\node[block] (rmp) {Solve RMP};

\node[block, right=of rmp] (heur)
{Heuristic\\labeling};

\node[decision, right=of heur] (heurcheck)
{Negative\\columns found?};

\node[block, right=of heurcheck] (exact)
{Exact\\labeling};

\node[decision, right=of exact] (exactcheck)
{Negative\\columns found?};

\node[block, right=of exactcheck] (term)
{Terminate};

\draw[line] (rmp) -- (heur);
\draw[line] (heur) -- (heurcheck);
\draw[line] (heurcheck) -- node[above, font=\scriptsize] {No} (exact);
\draw[line] (exact) -- (exactcheck);
\draw[line] (exactcheck) -- node[above, font=\scriptsize] {No} (term);

\draw[line]
(heurcheck.north)
-- ++(0,0.65)
-| node[pos=0.32, above, font=\scriptsize]
{Yes: add columns}
(rmp.north);

\draw[line]
(exactcheck.south)
-- ++(0,-0.65)
-| node[pos=0.30, below, font=\scriptsize]
{Yes: add columns}
(rmp.south);

\end{tikzpicture}

\caption{Original CG framework. Heuristic labeling is attempted first,
while exact labeling is invoked when no improving column is found.
CG terminates only when exact labeling finds no improving column.}
\label{fig:original_cg}
\end{figure}
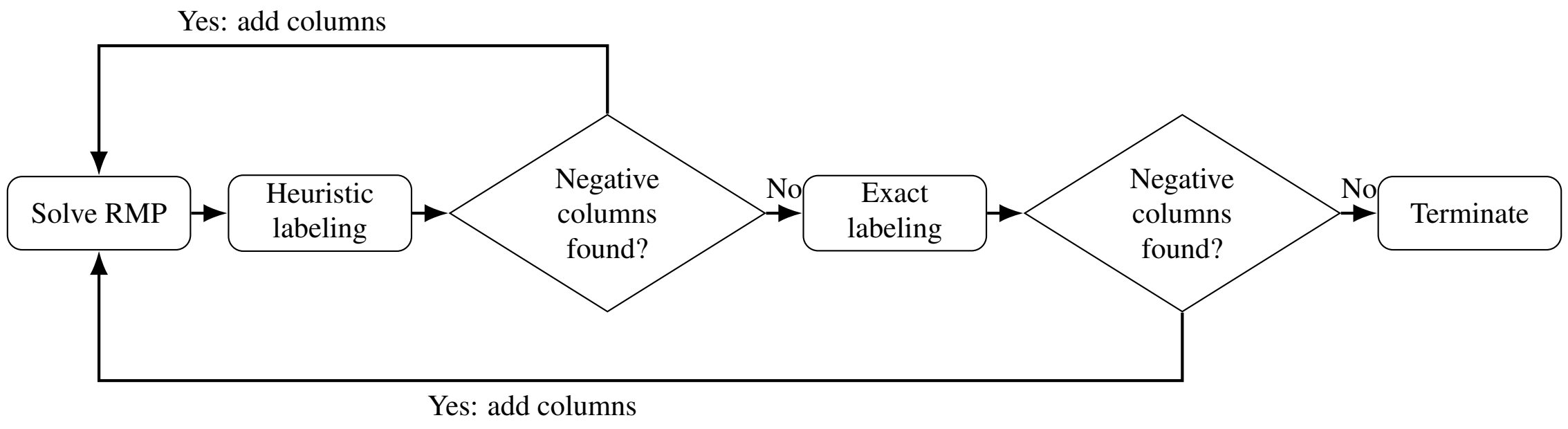

\subsection{Labeling Algorithm}\label{appendix:labeling_algorithm}

At CG iteration $t$, the RMP returns dual variables that define the pricing
problem $P_t$. A label $\ell$ represents a partial route and records its current
end node $i(\ell)$, visited customer set $V(\ell)$, accumulated resources, and
reduced cost $\bar c(\ell)$. Its feasible extension set is
\[
\mathcal F(\ell)
=
\left\{
j\in\mathcal N\setminus V(\ell):
\ell\oplus j \text{ is feasible}
\right\},
\]
where feasibility is determined by the capacity, time-window, and elementarity
constraints. Dominance removes labels that cannot lead to a better feasible
completion. The extension-ordering rule $g_t$ determines which active label is
processed next without changing the pricing formulation, feasible extensions,
reduced-cost calculation, or dominance rules.

\paragraph{Exact labeling.}
Exact labeling is the exact pricing procedure. It processes active
nondominated labels until the active label pool is empty and therefore certifies
whether any negative-reduced-cost column remains. The extension-ordering rule
may affect the number and order of labels explored, but not the optimal
reduced-cost value returned.

\begin{algorithm}[!htbp]
\caption{Exact labeling under extension ordering}
\label{alg:exhaustive_labeling}
\begin{algorithmic}[1]
\Require Dual variables and extension-ordering rule $g_t$
\Ensure Candidate column set $\mathcal C_t$

\State Initialize active label pool $\mathcal L_t$ and
$\mathcal C_t\leftarrow\emptyset$

\While{$\mathcal L_t\neq\emptyset$}
    \State Select and remove the highest-priority label $\ell$ under $g_t$
    \If{$\ell$ is a complete feasible route}
        \If{$\bar c(\ell)<0$}
            \State $\mathcal C_t\leftarrow\mathcal C_t\cup\{r(\ell)\}$
        \EndIf
    \Else
        \For{each $j\in\mathcal F(\ell)$}
            \State Generate $\ell'\leftarrow\ell\oplus j$
            \If{$\ell'$ is not dominated}
                \State Remove labels dominated by $\ell'$
                \State Add $\ell'$ to $\mathcal L_t$
            \EndIf
        \EndFor
    \EndIf
\EndWhile

\State \Return $\mathcal C_t$
\end{algorithmic}
\end{algorithm}

\paragraph{Early-stop labeling.}
Early-stop labeling is used as the heuristic pricing procedure in this paper. It follows the same extension, feasibility, reduced-cost, and dominance rules as exact labeling, but terminates once a prescribed stopping criterion
$\Gamma$ is satisfied, reaching a labeling budget. Hence, the extension-ordering rule can affect both the number of labels explored and the candidate columns discovered before termination.

\begin{algorithm}[!htbp]
\caption{Early-stop labeling under extension ordering}
\label{alg:early_stop_labeling}
\begin{algorithmic}[1]
\Require Dual variables, extension-ordering rule $g_t$, stopping criterion $\Gamma$
\Ensure Candidate column set $\mathcal C_t$

\State Initialize active label pool $\mathcal L_t$ and
$\mathcal C_t\leftarrow\emptyset$

\While{$\mathcal L_t\neq\emptyset$ \textbf{and} $\Gamma$ is not satisfied}
    \State Select and remove the highest-priority label $\ell$ under $g_t$
    \If{$\ell$ is a complete feasible route}
        \If{$\bar c(\ell)<0$}
            \State $\mathcal C_t\leftarrow\mathcal C_t\cup\{r(\ell)\}$
        \EndIf
    \Else
        \For{each $j\in\mathcal F(\ell)$}
            \State Generate $\ell'\leftarrow\ell\oplus j$
            \If{$\ell'$ is not dominated}
                \State Remove labels dominated by $\ell'$
                \State Add $\ell'$ to $\mathcal L_t$
            \EndIf
        \EndFor
    \EndIf
\EndWhile

\State \Return $\mathcal C_t$
\end{algorithmic}
\end{algorithm}

Other heuristic labeling strategies can also be combined with early stopping. In this paper, we focus on early-stop labeling so that heuristic and exact labeling share the same underlying labeling mechanism and differ only in their stopping condition.

\subsection{Resource-Based Priority Rules}
\label{appendix:priority_rules}

We consider five extension-ordering rules:
\[
\{
\text{minimum reduced cost},
\text{minimum load},
\text{minimum time},
\text{maximum load},
\text{maximum time}
\}.
\]
Each rule reflects a different search preference in the labeling process.

\begin{itemize}
    \item \textbf{Minimum reduced cost:}
    prioritizes partial routes that are currently most promising with respect to
    the pricing objective, with the goal of discovering negative-reduced-cost
    routes earlier.

    \item \textbf{Minimum load:}
    prioritizes labels with lower accumulated load, leaving more residual vehicle
    capacity and therefore more flexibility for future customer extensions.

    \item \textbf{Maximum load:}
    prioritizes labels with higher accumulated load, encouraging the search toward
    more resource-saturated and potentially near-complete routes.

    \item \textbf{Minimum time:}
    prioritizes labels with smaller accumulated time, preserving more remaining
    temporal flexibility for feasible future extensions.

    \item \textbf{Maximum time:}
    prioritizes labels with larger accumulated time, favoring deeper routes that
    are closer to exhausting their temporal resources or completing the route.
\end{itemize}

These rules modify only the order in which active labels are processed. Under exact labeling, they preserve pricing optimality while potentially changing
the labeling effort; under early-stop labeling, they can additionally affect the
set of improving columns discovered before termination.

\subsection{Motivation Experiments}
\label{appendix:motivation}

We conduct preliminary experiments on 672 Solomon CVRPTW instances~\citep{solomon1987algorithms} to examine the
effect of extension-ordering strategies. We compare five fixed and five random
strategies.\footnote{A fixed strategy uses the same extension-ordering rule
throughout all CG iterations of an instance, whereas a random strategy samples
an extension-ordering rule at each CG iteration.}

\begin{figure}[!htbp]
    \centering
    \begin{subfigure}{0.48\textwidth}
        \centering
        \includegraphics[width=\linewidth]{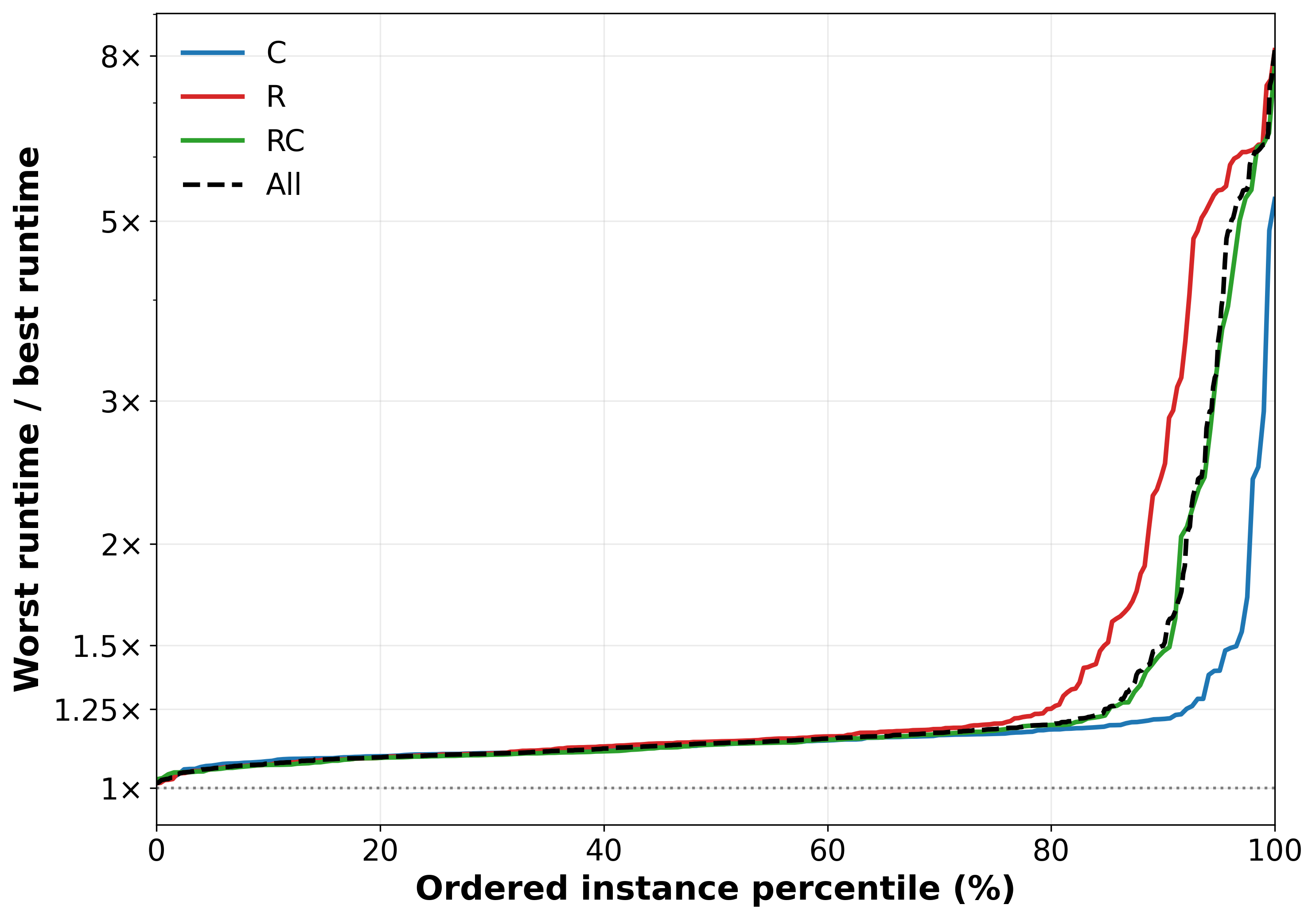}
        \caption{Ordered worst-to-best runtime ratio across extension-ordering strategies.}
        \label{fig:extension_strategy_effort}
    \end{subfigure}
    \hfill
    \begin{subfigure}{0.48\textwidth}
        \centering
        \includegraphics[width=\linewidth]{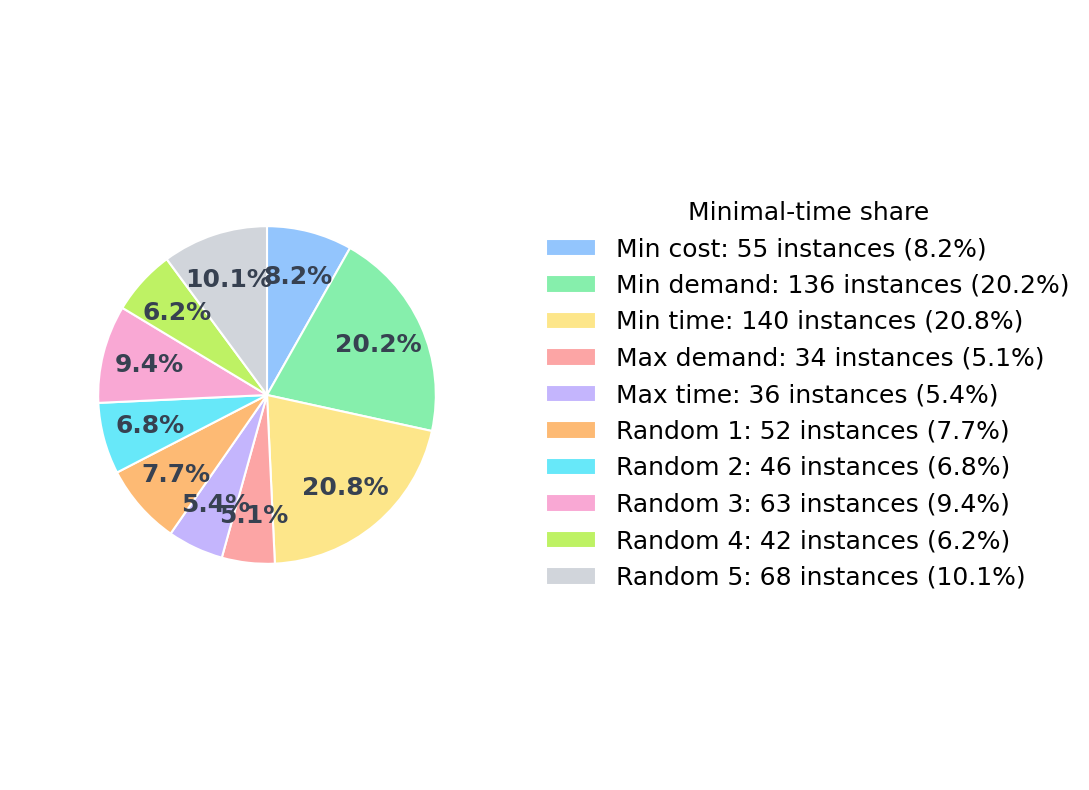}
        \caption{Fraction of instances for which each strategy achieves the minimum runtime.}
        \label{fig:extension_strategy_share}
    \end{subfigure}

    \vspace{1ex}

    \begin{subfigure}{0.78\textwidth}
        \centering
        \includegraphics[width=\linewidth]{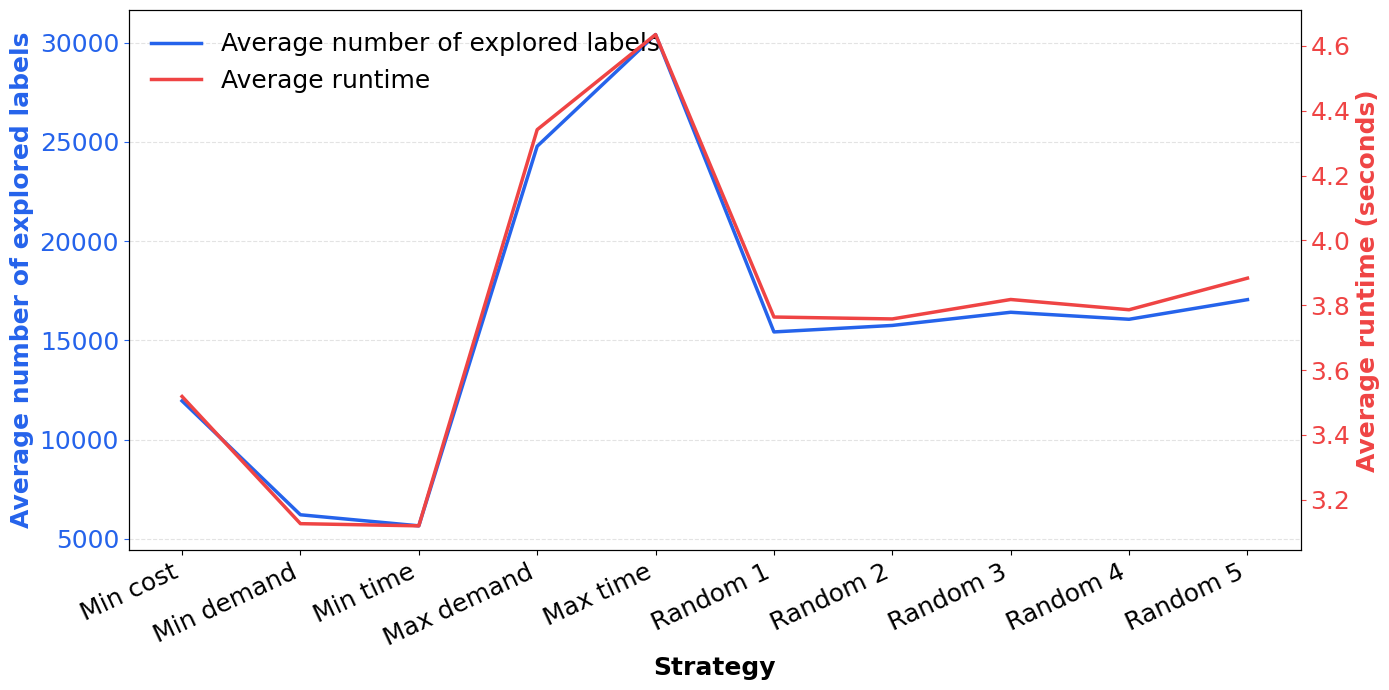}
        \caption{Average explored labels and average runtime under different extension-ordering strategies.}
        \label{fig:extension_strategy_surrogate}
    \end{subfigure}

    \caption{
    Empirical motivation for adaptive extension-ordering selection.
    Panel~(a) shows the instance-wise runtime variation across the ten candidate
    strategies, Panel~(b) shows that no single strategy dominates across
    instances, and Panel~(c) shows that the average number of explored labels
    closely tracks the average runtime across strategies. These results indicate
    both substantial potential gains from appropriate rule selection and support
    the use of explored-label count as a surrogate for pricing computational
    effort.
    }
    \label{fig:extension_strategy_motivation}
\end{figure}

Figure~\ref{fig:extension_strategy_motivation}(a) reports, for each instance, the
ratio between the largest and smallest runtimes obtained by the ten strategies,
ordered by percentile and shown separately for the \texttt{C}, \texttt{R}, and
\texttt{RC} classes as well as for all instances. The effect of priority-rule
selection varies substantially across instances. While the runtime differences
are small for many cases, nearly 20\% of the instances exhibit a worst-to-best
runtime ratio above 1.25, and the largest ratios approach eight. This indicates
that selecting an appropriate extension-ordering rule can have a substantial
impact on pricing efficiency.

Figure~\ref{fig:extension_strategy_motivation}(b) further shows that no single
strategy is consistently best. Each strategy achieves the minimum runtime on
some instances, while even the most frequently best strategy accounts for only
about 21\% of the instances. Together, the large within-instance runtime
variation and the absence of a uniformly dominant rule motivate selecting the
extension-ordering strategy adaptively according to the current CG state.

Figure~\ref{fig:extension_strategy_motivation}(c) compares, across the ten
strategies, the average number of explored labels and the average runtime.
The two curves exhibit closely aligned trends, indicating that the explored
label count tracks the actual computational effort of pricing well. This
supports using $N_t(g_t)$ as a practical surrogate for the pricing
computational cost $C(P_t)$.

%% file: Content/Methodology.tex
\section{Coupled-CG: Joint Learning for CG Acceleration}\label{sec:methodology}

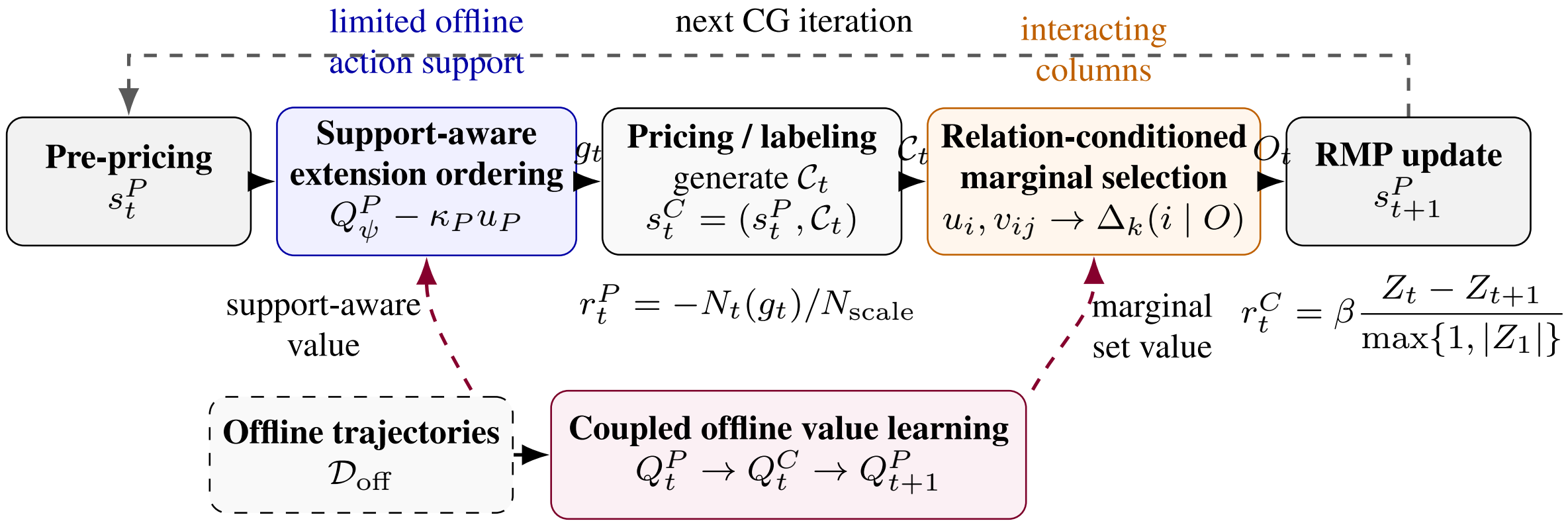
\begin{figure}[!htbp]
\centering

\resizebox{0.96\textwidth}{!}{%
\begin{tikzpicture}[
    >=Latex,
    every node/.style={
        font=\scriptsize,
        align=center
    },
    state/.style={
        rectangle,
        rounded corners,
        draw=black,
        fill=gray!10,
        minimum width=2.0cm,
        minimum height=1.05cm,
        inner sep=3pt
    },
    pricing/.style={
        rectangle,
        rounded corners,
        draw=blue!65!black,
        fill=blue!6,
        minimum width=2.45cm,
        minimum height=1.20cm,
        inner sep=3pt
    },
    process/.style={
        rectangle,
        rounded corners,
        draw=black,
        fill=gray!5,
        minimum width=2.45cm,
        minimum height=1.20cm,
        inner sep=3pt
    },
    column/.style={
        rectangle,
        rounded corners,
        draw=orange!75!black,
        fill=orange!7,
        minimum width=2.65cm,
        minimum height=1.20cm,
        inner sep=3pt
    },
    value/.style={
        rectangle,
        rounded corners,
        draw=purple!65!black,
        fill=purple!6,
        minimum width=3.5cm,
        minimum height=1.05cm,
        inner sep=4pt
    },
    data/.style={
        rectangle,
        rounded corners,
        draw=black,
        dashed,
        fill=gray!5,
        minimum width=2.2cm,
        minimum height=0.95cm,
        inner sep=3pt
    },
    arrow/.style={
        ->,
        thick
    },
    couple/.style={
        ->,
        thick,
        dashed,
        draw=purple!70!black
    },
    feedback/.style={
        ->,
        thick,
        dashed,
        draw=gray!70!black
    }
]

\node[state] (pre) {
    \textbf{Pre-pricing}\\
    $s_t^P$
};

\node[pricing, right=2mm of pre] (ordering) {
    \textbf{Support-aware}\\
    \textbf{extension ordering}\\
    $Q_\psi^P-\kappa_Pu_P$
};

\node[process, right=2mm of ordering] (pricingproc) {
    \textbf{Pricing / labeling}\\
    generate $\mathcal C_t$\\
    $s_t^C=(s_t^P,\mathcal C_t)$
};

\node[column, right=2mm of pricingproc] (selection) {
    \textbf{Relation-conditioned}\\
    \textbf{marginal selection}\\
    $u_i,v_{ij}\rightarrow\Delta_k(i\mid O)$
};

\node[state, right=2mm of selection] (rmp) {
    \textbf{RMP update}\\
    $s_{t+1}^P$
};

\draw[arrow] (pre) -- (ordering);

\draw[arrow]
    (ordering) --
    node[above, font=\scriptsize] {$g_t$}
    (pricingproc);

\draw[arrow]
    (pricingproc) --
    node[above, font=\scriptsize] {$\mathcal C_t$}
    (selection);

\draw[arrow]
    (selection) --
    node[above, font=\scriptsize] {$O_t$}
    (rmp);

\node[
    below=1mm of pricingproc,
    font=\scriptsize
] (rp) {
    $r_t^P=-N_t(g_t)/N_{\mathrm{scale}}$
};

\node[
    below=1mm of rmp,
    font=\scriptsize
] (rc) {
    $r_t^C=
    \beta\dfrac{Z_t-Z_{t+1}}
    {\max\{1,|Z_1|\}}$
};

\node[
    value,
    below=11mm of pricingproc,
    xshift=3mm
] (coupled) {
    \textbf{Coupled offline value learning}\\
    $Q_t^P\rightarrow Q_t^C\rightarrow Q_{t+1}^P$
};

\node[
    data,
    left=3mm of coupled
] (offline) {
    \textbf{Offline trajectories}\\
    $\mathcal D_{\mathrm{off}}$
};

\draw[arrow] (offline) -- (coupled);

\draw[couple]
    ([xshift=-26mm]coupled.north)
    to[out=120,in=-90]
    node[left, font=\scriptsize]
    {support-aware\\value}
    (ordering.south);

\draw[couple]
    ([xshift=20mm]coupled.north)
    to[out=60,in=-90]
    node[right, font=\scriptsize]
    {marginal\\set value}
    (selection.south);

\draw[feedback]
    (rmp.north)
    -- ++(0,5mm)
    -|
    node[pos=0.24, above=0.5mm, font=\scriptsize]
    {next CG iteration}
    (pre.north);

\node[
    above=0.8mm of ordering,
    font=\scriptsize,
    text=blue!65!black
] {
    limited offline\\
    action support
};

\node[
    above=0.8mm of selection,
    font=\scriptsize,
    text=orange!75!black
] {
    interacting\\
    columns
};

\end{tikzpicture}%
}

\caption{
Overview of \textsc{Coupled-CG}. Each CG iteration couples
extension ordering and column selection through
$Q_t^P \rightarrow Q_t^C \rightarrow Q_{t+1}^P$,
\rev{where $Q_t^P$ and $Q_t^C$ denote the long-term action-value
estimates for the pricing extension-ordering and column-selection
policies, respectively,} to accelerate CG.
}
\label{fig:coupled_cg_framework}

\end{figure}

Solving $\mathcal M$ from offline CG trajectories requires coordinating the two
decision stages, handling uneven support for extension-ordering actions, and
selecting interacting columns from a combinatorial candidate set.
\textsc{Coupled-CG} addresses these issues through coupled offline value
learning, support-aware extension ordering, and relation-conditioned marginal
column selection. Figure~\ref{fig:coupled_cg_framework} summarizes the overall
framework.

\subsection{Two-Stage MDP for \textsc{Coupled-CG} Decisions}
\label{sec:two_stage_mdp}

We formulate the coupled pricing and column-selection decisions as a two-stage MDP that directly follows the interaction identified in Section~\ref{sec:problem_setup}. At iteration $t$, the current RMP defines the pre-pricing state $s_t^P$, represented as a bipartite graph between column and constraint nodes with instance- and iteration-level features (Details are in Section~\ref{appendix:node_feature}). The pricing policy first selects an extension-ordering rule
\begin{equation}
g_t\sim\pi^P(\cdot\mid s_t^P),\qquad g_t\in\mathcal G,
\end{equation}
which determines the labeling search and candidate set $\mathcal C_t$, producing
\begin{equation}
s_t^C=\mathrm{PricingTransition}(s_t^P,g_t)
=(s_t^P,\mathcal C_t).
\end{equation}
The column-selection policy then chooses $O_t\subseteq\mathcal C_t$, and the selected columns update the RMP:
\begin{equation}
O_t\sim\pi^C(\cdot\mid s_t^C),
\qquad
s_{t+1}^P=\mathrm{RMPUpdate}(s_t^C,O_t).
\end{equation}

Thus, one CG iteration forms the sequential transition
$
s_t^P
\xrightarrow{\,g_t\,}
s_t^C
\xrightarrow{\,O_t\,}
s_{t+1}^P,
$
which explicitly couples pricing-search control with the downstream RMP update. A detailed specification of the states, actions, transitions,
rewards, and value functions for both decision stages is
provided in Section~\ref{app:mdp_specification}.

\paragraph{Stage-specific reward design.}
The two decisions affect CG differently: extension ordering
controls labeling effort, while column selection improves the
RMP and influences subsequent pricing problems. Accordingly,
we define stage-specific rewards as
\begin{equation}
r_t^P
=
-\frac{N_t(g_t)}{N_{\mathrm{scale}}},
\qquad
r_t^C
=
\beta
\frac{Z_t-Z_{t+1}}
{\max\{1,|Z_1|\}}.
\end{equation}

The pricing reward \rev{$r_t^P$} penalizes the number of labels
\rev{$N_t(g_t)$} explored \rev{at iteration $t$},
including heuristic labeling and exact fallback\rev{, so its sum over
an episode is the scaled cumulative labeling effort}.
Here, $N_{\mathrm{scale}}>0$ is a fixed normalization constant
estimated from offline training trajectories.
The column-selection reward \rev{$r_t^C$} measures the RMP objective
improvement, \rev{where $Z_t$ and $Z_{t+1}$ are the objective values
before and after the update.} \rev{It gives column selection feedback
at every iteration, whereas the labeling cost of a selected column
set appears only in later pricing calls.} The initial objective $Z_1$
provides normalization, and $\beta\geq0$ controls the reward
weight.

These rewards capture immediate computational effort and
optimization progress, while the coupled value functions
account for their long-term effects. \rev{Under} a
fixed initial RMP, undiscounted episodic returns, and exact
CG termination, the cumulative RMP-improvement reward
telescopes to \rev{$\beta(Z_1-Z_{T+1})/\max\{1,|Z_1|\}$, which is the
same for every policy because all policies reach the same LP
optimum $Z_{T+1}$}. Thus, it does
not alter the objective of minimizing cumulative labeling
effort $\mathcal M$, as formalized in Proposition~\ref{prop:reward_consistency}.

\begin{proposition}[Reward consistency with the acceleration objective]
\label{prop:reward_consistency}
For an undiscounted CG episode terminating under exact pricing, if
$N_{\mathrm{scale}}$ is fixed across policies, then
\begin{equation}
\max_{\pi}\sum_{t=1}^{T}(r_t^P+r_t^C)
\quad\Longleftrightarrow\quad
\min_{\pi}\sum_{t=1}^{T}N_t(g_t).
\end{equation}
Hence, the RMP-improvement reward provides dense shaping\rev{~\citep{ng1999policy}} without changing the optimizer of $\mathcal M$.
\end{proposition}

\noindent
The proof is provided in
Appendix~\ref{appendix:reward_consistency_proof}. Proposition~1 applies to the undiscounted episodic return before introducing the conservative learning regularizers, under a fixed initial RMP and exact CG termination. It shows that the proposed stage rewards are consistent with the cumulative labeling-effort objective. The conservative and support-aware terms introduced for offline learning, however, modify the value-learning procedure and are not covered by this equivalence result.

\subsection{Coupled Offline Value Learning}
\label{sec:coupled_value_learning}

We learn the two decision policies from pre-collected CG trajectories
generated by a randomized behavior policy. The details of trajectory construction are shown in Section~\ref{appendix:offline_dataset}. Each transition records
the pre-pricing state, the logged extension-ordering action, the
generated candidate columns, the logged selected column set, the
stage-specific rewards, and the next pre-pricing state. The dataset
is fixed before training, so policy learning does not require
additional interactions with the CG solver.

\rev{We define $Q_\psi^P(h_t,g)$ as the long-term action value of
extension-ordering rule $g$ and
$Q_\phi^C(\widetilde h_t,O)$ as the long-term value of selecting
column set $O$, where $h_t$ and $\widetilde h_t$ are the learned
representations of the pre-pricing and post-pricing states.} Because the decisions occur sequentially, their
Bellman targets are coupled across stages:
\begin{equation}
y_t^P
=
r_t^P
+(1-\mathrm{done}_t)
V_{\bar\phi}^C(\widetilde h_t),
\qquad
y_t^C
=
r_t^C
+(1-\mathrm{done}_t)
V_{\bar\psi}^P(h_{t+1}).
\end{equation}

\rev{The pricing target $y_t^P$ incorporates the downstream value of selecting
columns generated under the current extension-ordering rule.
The column-selection target $y_t^C$ incorporates the value of the next
pricing decision after the RMP update.} Thus, both decisions are
evaluated through their longer-term effects on the CG trajectory,
rather than through their immediate rewards alone.

The target values are constructed as
\begin{equation}
V_{\bar\phi}^C(\widetilde h_t)
=
Q_{\bar\phi,k}^C(\widetilde h_t,O_t^{\bar\phi}), \qquad 
V_{\bar\psi}^P(h_{t+1})=
\max_{g\in\mathcal G}
\left[
Q_{\bar\psi}^P(h_{t+1},g)
-\kappa_Pu^P(h_{t+1},g)
\right],
\end{equation}
where $O_t^{\bar\phi}$ is constructed by the target
column-selection policy. The terminal indicator suppresses
bootstrapping when the CG episode terminates.

Both value functions are trained using temporal-difference
objectives and stage-specific conservative regularization\rev{~\citep{kumar2020conservative}},
with alternating updates and slowly updated target networks.
The full optimization procedure is provided in Section~\ref{appendix:learning_details}.

\subsection{Stage-Specific Policy Models}
\label{sec:stage_specific_policies}
\paragraph{Support-aware extension ordering.}
We use five complementary extension-ordering rules to induce diverse labeling
search trajectories 
$$\mathcal G= \{
\text{minimum reduced cost},
\text{minimum load},
\text{minimum time},
\text{maximum load},
\text{maximum time}
\}.$$ 
The reduced-cost rule follows the pricing objective, while the load- and
time-based rules steer the search across different resource regimes. Detailed definitions of the five extension-ordering rules are provided in
Section~\ref{appendix:priority_rules}.

To motivate adaptive extension ordering, we further conduct experiments comparing
different fixed rules across instances. The results, reported in
Section~\ref{appendix:motivation}, show that the preferred rule varies across
instances, motivating the need to select the ordering strategy adaptively. In addition, to make adaptive extension ordering reliable under uneven offline action
coverage, we penalize rules that are weakly supported in similar states. Let $\xi_g$ denote the prototype of rule $g$ constructed from its logged state embeddings. We define
\begin{equation}
\widetilde{\omega}_P(h_t,g)
=
\exp\left(
-\frac{\|\bar h_t-\xi_g\|_2^2}{2\sigma^2}
\right),
\qquad
\omega_P(h_t,g)
=
\frac{\widetilde{\omega}_P(h_t,g)}
{\sum_{g'\in\mathcal G}\widetilde{\omega}_P(h_t,g')+\epsilon},
\end{equation}
where $\bar h_t=h_t/(\|h_t\|_2+\epsilon)$. We define the support uncertainty as $u_P(h_t,g)=1-\omega_P(h_t,g),$ with prototype construction detailed in
Section~\ref{appendix:support_kernel}. The extension-ordering rule is selected as
\begin{equation}
g_t
=
\arg\max_{g\in\mathcal G}
\left[
Q_\psi^P(h_t,g)-\kappa_Pu_P(h_t,g)
\right].
\end{equation}
where $\kappa_P\geq0$ penalizes weakly supported actions. This preserves adaptive pricing-search control while reducing extrapolation toward unsupported actions.

\paragraph{Relation-conditioned marginal column selection.}
To model column interactions without constructing an explicit subset-action
catalogue, we perform column selection through sequential marginal decisions. We first encode the candidate set $\mathcal C_t$ as a relation graph. An edge-aware encoder produces candidate representations $z_{t,i}$, from which
unary and symmetric pairwise heads predict individual values $u_i$ and compatibilities $v_{ij}=v_{ji}$. Implementation details of the relation graph, edge-aware encoder, and marginal selection procedure are provided in Section~\ref{appendix:marginal_column_selection}. For target cardinality $k>1$, we define the column-set value as
\begin{equation}
Q_{\phi,k}^C(\widetilde h_t,O)
=
\frac{1}{k}\sum_{i\in O}u_i
+
\frac{2}{k(k-1)}
\sum_{\substack{i<j\\ i,j\in O}}v_{ij}.
\end{equation}
For $k=1$, we use the unary value only. Under this factorization, the marginal
value of adding candidate $i\notin O$ is
\begin{equation}
\Delta_k(i\mid O)
=
\frac{1}{k}u_i
+
\frac{2}{k(k-1)}
\sum_{j\in O}v_{ij}.
\end{equation}
Starting from the mandatory minimum-reduced-cost column, the policy repeatedly selects $i^\star$ using $\arg\max_{i\in\mathcal C_t\setminus O}
\Delta_k(i\mid O)$ until $k$ columns are selected or no unused candidate remains. The marginal construction avoids enumerating $\binom{M}{k}$ subsets while capturing column quality and pairwise interactions, and naturally generalizes across different $M$ and $k$.

\paragraph{Training and inference.}
The two value functions are trained with temporal-difference and stage-specific conservative losses\rev{~\citep{kumar2020conservative}},
\begin{equation}
\mathcal L^P
=
\mathcal L_{\mathrm{TD}}^P
+
\lambda_P\mathcal L_{\mathrm{CCQL}}^P,
\qquad
\mathcal L^C
=
\mathcal L_{\mathrm{TD}}^C
+
\lambda_C\mathcal L_{\mathrm{MCR}}^C,
\end{equation}
where $\mathcal L_{\mathrm{CCQL}}^P$ and $\mathcal L_{\mathrm{MCR}}^C$ denote the conservative pricing and marginal-selection regularizers, respectively. Together, conservative support-aware pricing, relation-conditioned marginal selection, and coupled Bellman updates provide robust and scalable coordination of both decisions toward the common objective $\mathcal M$; optimization details are given in Section~\ref{appendix:learning_details}.

\subsection{Bipartite Graph Representation} \label{appendix:node_feature}
Node features used for the capacitated vehicle routing problem with time windows (CVRPTW):
\begin{enumerate}
\item \textbf{Column node features $\mathcal X$:}
\rev{Following the column feature representation introduced by
\citet{chi2022deep}, we characterize each route variable using
its optimization status, structural properties, and historical
basis information. These features capture both the current RMP
solution and the evolution of individual columns across CG
iterations, with adaptations to the vehicle routing setting.}

\begin{enumerate}
    \item \textbf{Reduced cost:}
    Reduced cost is a quantity measuring the net benefit of
    adding a route to the current solution: it compares the
    route’s actual operating cost with the value gained from
    serving customers that are currently expensive to cover,
    while also accounting for the opportunity cost of using
    one vehicle. A route has a negative reduced cost only if
    the benefit of covering those customers outweighs both
    its travel cost and the scarcity of available vehicles.

    \item \textbf{Connectivity of column node:}
    \rev{The degree of a column node in the bipartite graph,
    measured by the number of constraints in which the
    corresponding route has a nonzero coefficient. For
    customer-coverage constraints, this feature reflects
    the number of customers served by the route.}

    \item \textbf{Solution value:}
    \rev{The current primal value of the route variable obtained
    from RMP optimization. For candidate routes that have
    not yet been incorporated into the RMP, this feature
    is initialized to zero.}

    \item \textbf{Route cost:}
    The route cost represents the total operational cost
    incurred by executing a specific route, typically
    including travel distance or time as well as any
    penalties for violating time windows, waiting, or other
    service-related considerations. In the CVRPTW context,
    it reflects the true resource consumption of assigning
    one vehicle to serve the sequence of customers defined
    by that route.

    \item \rev{\textbf{Number of iterations in the basis:}
    A historical counter measuring the duration of a
    column's membership in the RMP basis. This feature
    captures the persistence of its basic status across
    successive CG iterations.}

    \item \rev{\textbf{Number of iterations outside the basis:}
    A historical counter measuring the duration for which
    a column remains nonbasic, providing information about
    its previous participation in the RMP basis.}

    \item \rev{\textbf{Basis exit indicator:}
    A binary indicator equal to one if the column
    transitioned from basic to nonbasic status during
    the most recent RMP update, and zero otherwise.}

    \item \rev{\textbf{Basis entry indicator:}
    A binary indicator equal to one if the column
    transitioned from nonbasic to basic status during
    the most recent RMP update, and zero otherwise.}
\end{enumerate}

\item \textbf{Constraint node features $\mathcal Y$:}
\rev{Each constraint node represents a row of the RMP.
Since CG introduces variables without modifying the
constraint set, the number of constraint nodes remains
unchanged throughout the CG procedure for a given instance.
We adopt the dual value and connectivity features from
\citet{chi2022deep} and augment them with instance-specific
customer attributes.}

\begin{enumerate}
    \item \textbf{Dual value:}
    \rev{The optimal dual multiplier associated with the
    corresponding RMP constraint at the current iteration.
    This feature captures the marginal value of the
    constraint and contributes to the reduced-cost
    calculation in the subsequent pricing problem.}

    \item \textbf{Connectivity of constraint node:}
    \rev{The degree of a constraint node, computed as the
    number of column variables having nonzero coefficients
    in the corresponding RMP row.}

    \item \textbf{Instance information:}
    Customer-related attributes including demand, service
    duration, time window bounds (earliest and latest
    allowable service times), and spatial location
    represented by node coordinates $(x,y)$.
\end{enumerate}
\end{enumerate}
In addition to node features, the bipartite graph representation also incorporates edge features $\mathcal E$ to characterize the relationships between column nodes and constraint nodes:
\begin{itemize}
    \item \textbf{Incidence}: Always equal to $1.0$ for an existing edge. This feature simply indicates that the corresponding column participates in the associated constraint.
    
    \item \textbf{Constraint type}: Encoded as $1.0$ if the edge connects to the convexity constraint, and $0.0$ otherwise. This feature allows the model to distinguish between customer-coverage constraints and the convexity constraint.
    
    \item \textbf{Customer ID}: The index of the customer associated with a customer-coverage constraint. For edges connected to non-customer constraints, such as convexity constraints, this feature is set to $0$.
\end{itemize}

\subsection{Two-Stage MDP Specification}
\label{app:mdp_specification}

Table~\ref{tab:two_stage_mdp} summarizes the two-stage Markov
decision process (MDP) underlying \textsc{Coupled-CG}. Each CG
iteration consists of a pricing decision followed by a
column-selection decision. The pricing stage determines the
labeling search effort and the candidate columns available for
selection, while the column-selection stage updates the RMP and
influences the subsequent pricing problem.

\begin{table}[!htbp]
\centering
\caption{Summary of the two-stage MDP in \textsc{Coupled-CG}.}
\label{tab:two_stage_mdp}
\small
\begin{tabular}{p{0.16\linewidth} p{0.38\linewidth} p{0.38\linewidth}}
\toprule
\textbf{MDP element}
& \textbf{Pricing stage ($P$)}
& \textbf{Column-selection stage ($C$)} \\
\midrule

State
& $s_t^P$: current RMP state, including existing columns,
dual variables, and routing-instance characteristics.
& $s_t^C=(s_t^P,\mathcal C_t)$: post-pricing state,
including the current RMP and candidate columns generated
by pricing. \\

State representation
& $h_t=\operatorname{GNN}_{\theta}(s_t^P)$:
graph-based representation of the current RMP used
to guide extension ordering.
& $\widetilde h_t=\operatorname{Enc}_{\theta}^{C}
(h_t,\mathcal C_t)$:
representation incorporating candidate-column features
and their relationships for column selection. \\

Action
& $g_t\in\mathcal G$:
select one of five label-extension ordering rules
to guide the labeling search.
& $O_t\subseteq\mathcal C_t$, $|O_t|\leq k$:
select up to $k$ candidate columns
for insertion into the RMP. \\

Transition
& $s_t^C=\operatorname{PricingTransition}(s_t^P,g_t)$:
execute labeling under the selected rule to generate
the candidate-column set $\mathcal C_t$.
& $s_{t+1}^P=\operatorname{RMPUpdate}(s_t^C,O_t)$:
insert selected columns and reoptimize the RMP
to obtain the next pricing state. \\

Immediate effect
& Determines labeling effort $N_t(g_t)$
and the composition of the generated candidate
set $\mathcal C_t$.
& Determines RMP objective improvement
$Z_t-Z_{t+1}$ and updated dual variables
for subsequent pricing. \\

Reward
& $r_t^P=-N_t(g_t)/N_{\mathrm{scale}}$:
penalizes normalized labeling effort.
& $r_t^C=\beta(Z_t-Z_{t+1})/
\max\{1,|Z_1|\}$:
measures normalized RMP objective improvement. \\

Action-value function
& $Q_\psi^P(h_t,g_t)$:
estimates the long-term value of an
extension-ordering decision.
& $Q_{\phi,k}^C(\widetilde h_t,O_t)$:
estimates the long-term value of a selected
column subset, accounting for column interactions. \\

Next stage
& Proceeds to column selection at $s_t^C$
when valid candidate columns are available.
& Proceeds to pricing at $s_{t+1}^P$
after the RMP update, unless CG terminates. \\

\bottomrule
\end{tabular}
\end{table}

The two stages form the sequential transition
\begin{equation}
s_t^P \xrightarrow{g_t} s_t^C
\xrightarrow{O_t} s_{t+1}^P.
\end{equation}
The pricing reward directly penalizes computational effort,
whereas the column-selection reward provides intermediate
feedback on RMP progress. Their longer-term effects are
captured through the coupled Bellman targets introduced
in Section~4.2. CG terminates only when exact pricing
certifies that no negative-reduced-cost column remains.

\subsection{Offline Dataset Construction}
\label{appendix:offline_dataset}

We construct the offline dataset by running randomized column-generation
trajectories on sampled instances from the C family of the Solomon CVRPTW
benchmark. The resulting dataset contains $300$ instances. Each instance forms one complete episode with a
terminal record. We use a candidate-pool size of $M=10$ and select at most
$k=5$ columns per iteration.

\paragraph{Behavior trajectory generation.}
For each training instance, we run a randomized CG behavior policy and record
one transition per CG iteration. At iteration $t$, the RMP is first solved and
its current state $s_t$ is recorded. An extension-ordering rule
$\widehat g_t\in\mathcal G$ is then sampled and used for pricing. Heuristic
labeling is applied first, with exact labeling used when no improving route is
found. The resulting novel negative-reduced-cost routes form the candidate set
$\widehat{\mathcal C}_t$.

The candidate routes are ordered by reduced cost. The most negative-reduced-cost route is retained as the first selected column,
and the remaining columns are sampled sequentially without replacement,
yielding the logged selection sequence
$\widehat{\boldsymbol O}_t
=
(\widehat i_{t,1},\ldots,\widehat i_{t,k})$.
The corresponding selected set
$\widehat O_t
=
\{\widehat i_{t,1},\ldots,\widehat i_{t,k}\}$
is then added to the RMP, which is reoptimized to obtain
$s_{t+1}$. If pricing produces no improving route or no column can be added,
the transition is marked terminal.

\paragraph{Logged transition.}
Each CG iteration produces an offline raw transition 
\begin{equation}
\left(
s_t,\widehat g_t,
\widehat{\mathcal C}_t,\widehat O_t,
m_t,s_{t+1},\mathrm{done}_t
\right),
\end{equation}
where $\widehat g_t$ is the logged extension-ordering action,
$\widehat{\mathcal C}_t$ is the generated candidate set,
$\widehat O_t$ is the logged selected subset, and
$\mathrm{done}_t$ indicates whether the CG trajectory terminates. 

The raw metric vector is
\begin{equation}
\label{eq:raw_metrics}
m_t
=
\left(
N_t^{\mathrm{pre}},
N_t^{\mathrm{cur}},
Z_t^{\mathrm{pre}},
Z_t^{\mathrm{cur}}
\right),
\end{equation}
where $N_t^{\mathrm{pre}}$ and $N_t^{\mathrm{cur}}$ denote the recorded
label-extension counts before and during the current pricing update, and
$Z_t^{\mathrm{pre}}$ and $Z_t^{\mathrm{cur}}$ denote the RMP objective values
before and after the corresponding RMP update. These measurements are used to
construct the stage-specific rewards during offline training.

During dataset aggregation, an entire source instance is discarded if either
$Z_t^{\mathrm{pre}}$ or $Z_t^{\mathrm{cur}}$ contains a missing, non-numeric,
NaN, or infinite value. 

\paragraph{Reward construction.}
For the pricing stage, we use the current labeling effort
$N_t^{\mathrm{cur}}$ and define
\begin{equation}
r_t^P
=
-
\frac{N_t^{\mathrm{cur}}}
{N_{\mathrm{scale}}},
\end{equation}
where
\begin{equation}
N_{\mathrm{scale}}
=
\frac{1}{|\mathcal I_{\mathrm{off}}|}
\sum_{\mathcal I\in\mathcal I_{\mathrm{off}}}
N(\mathcal I),
\end{equation}
$\mathcal I_{\mathrm{off}}$ denotes the set of offline training instances, and
$N(\mathcal I)$ is the total labeling effort over the complete CG trajectory
of instance $\mathcal I$. Thus, $N_{\mathrm{scale}}$ is the average total
labeling effort across the offline episodes.

For the column-selection stage, we use the RMP objective improvement induced by
the selected columns:
\begin{equation}
r_t^C
=
\beta
\frac{
Z_t^{\mathrm{pre}}-Z_t^{\mathrm{cur}}
}{
\max\{1,|Z_1|\}
},
\end{equation}
where $Z_1$ is the initial RMP objective and $\beta\geq0$ scales the
objective-improvement reward.

\paragraph{State and candidate features.}
The RMP state is represented as a bipartite graph between column and constraint
nodes. Column-node features include reduced cost, solution value, route cost,
connectivity, and basis-history information, while constraint-node features
include dual value, connectivity, customer demand, time-window information,
service time, and spatial coordinates. Edge features describe RMP incidence
and constraint type.

Each candidate route is represented by its reduced cost, route cost, number of
customers, total demand, route duration, normalized candidate rank, and overlap
with the current RMP. For relation-aware column selection, we additionally
construct a complete candidate-relation graph using pairwise customer overlap,
arc overlap, cosine distance, and Jaccard distance.

\subsection{Prototype-Based Local Action Support}
\label{appendix:support_kernel}

For each extension-ordering rule $g\in\mathcal G$, we maintain an
action-specific prototype $\xi_g$ in the learned state-embedding space. The
prototype summarizes the offline states for which rule $g$ is observed. Given
the normalized embedding
\begin{equation}
\bar h_j
=
\frac{h_j}{\|h_j\|_2+\epsilon},
\end{equation}
the prototype is initialized from the normalized centroid of the corresponding
logged states and updated during training by an exponential moving average. For a
minibatch $\mathcal B$, let
\begin{equation}
\mathcal B_g
=
\left\{
j\in\mathcal B:\widehat g_j=g
\right\},
\end{equation}
and define
\begin{equation}
\mu_g^{\mathcal B}
=
\frac{
\frac{1}{|\mathcal B_g|}
\sum_{j\in\mathcal B_g}\bar h_j
}{
\left\|
\frac{1}{|\mathcal B_g|}
\sum_{j\in\mathcal B_g}\bar h_j
\right\|_2+\epsilon
}.
\end{equation}
When $\mathcal B_g\neq\emptyset$, the prototype is updated as
\begin{equation}
\widetilde\xi_g
=
(1-\rho)\xi_g+\rho\mu_g^{\mathcal B},
\qquad
\xi_g
\leftarrow
\frac{\widetilde\xi_g}
{\|\widetilde\xi_g\|_2+\epsilon},
\end{equation}
where $\rho\in(0,1]$ is the prototype update rate.

For a current state representation $h_t$, we normalize
\begin{equation}
\bar h_t
=
\frac{h_t}{\|h_t\|_2+\epsilon}
\end{equation}
and compute its similarity to each action prototype using
\begin{equation}
\widetilde{\omega}_P(h_t,g)
=
\exp\left(
-\frac{\|\bar h_t-\xi_g\|_2^2}{2\sigma^2}
\right)
\mathbf 1\{n_g>0\},
\end{equation}
where $\sigma>0$ controls the similarity bandwidth and $n_g$ denotes the
number of logged training samples associated with rule $g$. The similarities
are normalized across extension-ordering rules:
\begin{equation}
\omega_P(h_t,g)
=
\frac{
\widetilde{\omega}_P(h_t,g)
}{
\sum_{g'\in\mathcal G}
\widetilde{\omega}_P(h_t,g')
+\epsilon
}.
\end{equation}
A larger $\omega_P(h_t,g)$ indicates stronger relative support for rule $g$
around the current state representation. We therefore define
\begin{equation}
u_P(h_t,g)
=
1-\omega_P(h_t,g)
\end{equation}
as the corresponding support uncertainty used in both conservative training
and support-adjusted action selection.

\subsection{Relation-Conditioned Marginal Column Selection}
\label{appendix:marginal_column_selection}

This section provides implementation details for the relation-conditioned
marginal column-selection policy introduced in the main text.

\paragraph{Candidate relation graph.}
For each pricing call, the generated candidate set
$\mathcal C_t=\{c_1,\ldots,c_M\}$ is represented as a complete relation graph.
Each candidate $c_i$ is initialized with its route-level features, including
reduced cost, route cost, number of customers, total demand, route duration,
normalized rank, and overlap with the current RMP. For each pair
$(c_i,c_j)$, we construct an edge feature vector
\begin{equation}
e_{ij}
=
\left[
o_{ij}^{\mathrm{cust}},
o_{ij}^{\mathrm{arc}},
d_{ij}^{\mathrm{cos}},
d_{ij}^{\mathrm{Jac}}
\right],
\end{equation}
where the four components denote customer overlap, arc overlap, cosine
distance, and Jaccard distance, respectively. These relations provide
information about redundancy and complementarity among generated columns.

\paragraph{Relation-aware encoding and value heads.}
An edge-aware graph encoder propagates information over the candidate-relation
graph and produces a representation $z_{t,i}$ for each candidate. The
pre-pricing representation $h_t$ is incorporated through
$\operatorname{Enc}_{\theta}^{C}$, yielding the post-pricing representation
$\widetilde h_t$. A unary head assigns an individual value
\begin{equation}
u_i=f_u(z_{t,i},\widetilde h_t),
\end{equation}
while a symmetric pairwise head evaluates the interaction between candidates,
\begin{equation}
v_{ij}
=
\frac{1}{2}
\left[
f_v(z_{t,i},z_{t,j},e_{ij})
+
f_v(z_{t,j},z_{t,i},e_{ji})
\right],
\end{equation}
so that $v_{ij}=v_{ji}$.

\paragraph{Marginal selection.}
Selection starts from the mandatory minimum reduced cost column. At each
subsequent step, the marginal value $\Delta_k(i\mid O)$ defined in the main
text is evaluated for every unused candidate, and the candidate with the
largest marginal value is added to $O$. The procedure continues until $k$
columns have been selected or no unused candidate remains.

Unlike an explicit subset-action representation, which requires evaluating up
to $\binom{M}{k}$ candidate subsets, marginal selection evaluates individual
remaining candidates sequentially. Given precomputed unary and pairwise
values, selecting $k$ columns requires at most $O(kM)$ marginal comparisons.

\paragraph{Variable candidate-pool and selection sizes.}
The policy does not associate actions with predefined subset indices.
Instead, unary and pairwise values are computed directly from the candidates
present in the current pricing call. Consequently, the same network
parameterization can be applied when the candidate-pool size $M$ or target
selection cardinality $k$ changes, without redefining the output action
space.

\subsection{Learning Details}
\label{appendix:learning_details}

\paragraph{Temporal-difference objectives.}
The two value functions are trained using the coupled Bellman targets defined in
Section~4.2. For the extension-ordering stage, the temporal-difference loss is
\begin{equation}
\mathcal L_{\mathrm{TD}}^P
=
\mathbb E_{\mathcal D_{\mathrm{off}}}
\left[
\left(
Q_\psi^P(h_t,\widehat g_t)-y_t^P
\right)^2
\right].
\end{equation}
For the column-selection stage, the logged selected set is evaluated using the
relation-conditioned set value:
\begin{equation}
\mathcal L_{\mathrm{TD}}^C
=
\mathbb E_{\mathcal D_{\mathrm{off}}}
\left[
\left(
Q_{\phi,k}^C(\widetilde h_t,\widehat O_t)-y_t^C
\right)^2
\right],
\end{equation}
where
\begin{equation}
k
=
\min\{k,|\widehat{\mathcal C}_t|\}.
\end{equation}

\paragraph{Conservative regularization for extension ordering.}
For the extension-ordering stage, unreliable counterfactual values are
regularized using
\begin{equation}
\mathcal L_{\mathrm{CCQL}}^P
=
\mathbb E_{\mathcal D_{\mathrm{off}}}
\left[
\log
\sum_{g\in\mathcal G}
\exp
\left(
Q_{\psi}^P(h_t,g)
+
\kappa_P u_P(h_t,g)
\right)
-
Q_{\psi}^P(h_t,\widehat g_t)
\right].
\end{equation}
The uncertainty term $u_P(h_t,g)$ places stronger conservative pressure on
high values assigned to weakly supported actions. The extension-ordering loss is
\begin{equation}
\mathcal L^P
=
\mathcal L_{\mathrm{TD}}^P +
\lambda_P\mathcal L_{\mathrm{CCQL}}^P,
\end{equation}
where $\lambda_P\geq0$ controls the strength of conservative regularization.

\paragraph{Conservative regularization for marginal column selection.}
The column-selection policy constructs a set through sequential marginal
decisions rather than treating all subsets as explicit actions. We therefore apply conservative regularization to the same marginal decisions
used by the policy. For the logged selection sequence
$\widehat{\boldsymbol O}_t$ defined in
Section~\ref{appendix:offline_dataset}, let
\begin{equation}
\widehat O_{t,\ell}
=
\left\{
\widehat i_{t,1},
\ldots,
\widehat i_{t,\ell}
\right\}
\end{equation}
denote the set of columns selected within the first $\ell$ marginal decisions. For $k>1$, the conservative loss is
\begin{equation}
\mathcal L_{\mathrm{MCR}}^C
=
\mathbb E_{\mathcal D_{\mathrm{off}}}
\left[
\frac{1}{k-1}
\sum_{\ell=1}^{k-1}
\left(
\log
\sum_{i\in
\widehat{\mathcal C}_t
\setminus
\widehat O_{t,\ell}}
\exp
\left(
\Delta_{k}
(i\mid\widehat O_{t,\ell})
\right)
-
\Delta_{k}
\left(
\widehat i_{t,\ell+1}
\mid
\widehat O_{t,\ell}
\right)
\right)
\right].
\end{equation}
For $k=1$, we set
$\mathcal L_{\mathrm{MCR}}^C=0$.
This regularization suppresses unsupported high marginal values while
preserving the values of logged marginal choices, without enumerating the
$\binom{M}{k}$ possible subsets.

The column-selection loss is
\begin{equation}
\mathcal L^C
=
\mathcal L_{\mathrm{TD}}^C
+
\lambda_C\mathcal L_{\mathrm{MCR}}^C,
\end{equation}
where $\lambda_C\geq0$ controls the strength of conservative regularization for
marginal column selection.

\paragraph{Coupled alternating training.}
Because $y_t^P$ depends on the target column-selection network and $y_t^C$
depends on the target extension-ordering network, the two value functions are
trained within one coupled offline-learning procedure. For each minibatch, the
target column-selection policy constructs
\begin{equation}
O_t^{\bar\phi}
=
\operatorname{MarginalSelect}
\left(
\widetilde h_t,
\widehat{\mathcal C}_t,
k;
\bar\phi
\right),
\end{equation}
using the marginal rule in Section~4.3. Its target value is
\begin{equation}
V_{\bar\phi}^C(\widetilde h_t)
=
Q_{\bar\phi,k}^C
\left(
\widetilde h_t,
O_t^{\bar\phi}
\right).
\end{equation}
The extension-ordering target
is
\begin{equation}
y_t^P
=
r_t^P
+(1-\text{done}_t)
V_{\bar\phi}^C(\widetilde h_t).
\end{equation}

After the selected columns update the RMP, the target extension-ordering value
for the next iteration is
\begin{equation}
V_{\bar\psi}^P(h_{t+1})
=
\max_{g\in\mathcal G}
\left[
Q_{\bar\psi}^P(h_{t+1},g)
-
\kappa_P u_P(h_{t+1},g)
\right],
\end{equation}
and the column-selection target is
\begin{equation}
y_t^C
=
r_t^C
+(1-\mathrm{done}_t)
V_{\bar\psi}^P(h_{t+1}).
\end{equation}
The extension-ordering and column-selection networks are updated alternately
using $\mathcal L^P$ and $\mathcal L^C$, respectively.

The target-network parameters are updated by Polyak averaging:
\begin{equation}
\bar\psi
\leftarrow
\tau_P\psi+(1-\tau_P)\bar\psi,
\qquad
\bar\phi
\leftarrow
\tau_C\phi+(1-\tau_C)\bar\phi,
\end{equation}
where $\tau_P,\tau_C\in(0,1]$ denote the corresponding target-update rates.
The shared state encoder is updated through both stage-specific losses during
alternating training.

The complete training and inference procedure is summarized in
Algorithm~\ref{alg:coupled_cql}.
\begin{algorithm}[!htbp]
\caption{Coupled Two-Stage Offline Conservative Q-Learning}
\label{alg:coupled_cql}
\begin{algorithmic}[1]

\Require Offline dataset $\mathcal D_{\mathrm{off}}$, strategy set
$\mathcal G$, target selection size $k$.
\Ensure Trained extension-ordering and column-selection policies.

\Statex \textbf{Offline training}
\State Initialize the shared encoder, $Q_\psi^P$, $Q_\phi^C$, and target networks.

\For{each training iteration}
    \State Sample minibatch $\mathcal B\subseteq\mathcal D_{\mathrm{off}}$
    and construct $r_t^P,r_t^C$ from $m_t$.

    \State Encode
    $h_t=\operatorname{GNN}_\theta(s_t^P)$ and
    $\widetilde h_t=
    \operatorname{Enc}_\theta^C(h_t,\widehat{\mathcal C}_t)$.

    \State Compute $u_P(h_t,g)$ for all $g\in\mathcal G$ and set
    $k=\min\{k,|\widehat{\mathcal C}_t|\}$.

    \State Construct
    $O_t^{\bar\phi}
    =
    \operatorname{MarginalSelect}
    (\widetilde h_t,\widehat{\mathcal C}_t,k;\bar\phi)$
    when $\widehat{\mathcal C}_t\neq\emptyset$.

    \State Compute the coupled pricing target $y_t^P$ and update
    $\theta,\psi$ by minimizing
    $\mathcal L^P
    =
    \mathcal L_{\mathrm{TD}}^P
    +
    \lambda_P\mathcal L_{\mathrm{CCQL}}^P$.

    \If{$\widehat{\mathcal C}_t\neq\emptyset$}
        \State Encode
        $h_{t+1}=\operatorname{GNN}_\theta(s_{t+1}^P)$.

        \State Compute
        $
        V_{\bar\psi}^P(h_{t+1})
        =
        \max_{g\in\mathcal G}
        [
        Q_{\bar\psi}^P(h_{t+1},g)
        -
        \kappa_Pu_P(h_{t+1},g)
        ].
        $

        \State Compute the coupled column-selection target $y_t^C$ and update
        $\theta,\phi$ by minimizing
        $\mathcal L^C
        =
        \mathcal L_{\mathrm{TD}}^C
        +
        \lambda_C\mathcal L_{\mathrm{MCR}}^C$.
    \EndIf

    \State Update target networks by Polyak averaging.
\EndFor

\Statex \textbf{Online inference}
\State Initialize the RMP with feasible columns.

\While{true}
    \State Solve the RMP and compute
    $h_t=\operatorname{GNN}_\theta(s_t^P)$.

    \State Select
    $
    g_t
    =
    \arg\max_{g\in\mathcal G}
    [
    Q_\psi^P(h_t,g)-\kappa_Pu_P(h_t,g)
    ].
    $

    \State Run heuristic labeling under $g_t$ to obtain $\mathcal C_t$;
    invoke exact labeling if $\mathcal C_t=\emptyset$.

    \If{$\mathcal C_t=\emptyset$}
        \State \textbf{break}
    \EndIf

    \State Compute
    $\widetilde h_t=
    \operatorname{Enc}_\theta^C(h_t,\mathcal C_t)$ and
    $k=\min\{k,|\mathcal C_t|\}$.

    \State
    $
    O_t
    =
    \operatorname{MarginalSelect}
    (\widetilde h_t,\mathcal C_t,k;\phi).
    $

    \State Add $O_t$ to the RMP.
\EndWhile

\State \Return Final RMP solution.

\end{algorithmic}
\end{algorithm}

%% file: Content/Experiment.tex
\section{Experiments}
We evaluate \textsc{Coupled CG} on Solomon CVRPTW instances from the six families C1, C2, R1, R2, RC1, and RC2~\citep{solomon1987algorithms}. All learning-based methods are trained on the same pool of 300 C-family instances with 6--18 customers, ensuring identical instance-level training access across methods. The competing learning-based methods construct their training data from the same instance pool according to their respective learning procedures. Details of dataset construction and baseline adaptation are provided in Sections~\ref{appendix:offline_dataset} and~\ref{appendix:baseline_dataset_preparation}.

\rev{We compare \textsc{Coupled-CG} with four baselines that share its CG backbone. Plain CG~\citep{desaulniers2006column} uses the minimum-reduced-cost extension order and adds the $k$ most-negative-reduced-cost columns. The other three learn column selection only and leave the extension order unchanged: Yuan CG~\citep{Yuan2024} by on-policy reinforcement learning, Morabit CG~\citep{Morabit2021} by a supervised classifier trained on expert labels, and FFCG~\citep{hu2025ffcg} by deep Q-learning. Sections~\ref{appendix:training_time}--\ref{appendix:neural_solver_comparison} add training times and comparisons with learning-based pricing~\citep{morabit2023,Xu2023RLHH,Abouelrous2025RLPricing}, adaptive stabilization~\citep{Shen2024AdaptiveCG}, and neural routing solvers~\citep{liu2024multitask,zhou2024mvmoe,berto2024routefinder}.}

Evaluation is conducted entirely on held-out instances. We consider two size regimes: small instances with 6--18 customers, which lie within the training size range, and larger instances with 20--60 customers, which lie outside the training size range and evaluate generalization to larger problems. Since training uses only C-family instances, evaluation on R and RC families additionally assesses transfer across instance structures. Performance is measured by CG time, end-to-end wall-clock time, the number of labels explored, the number of CG iterations, and the percentage of instances completed within the time limit. Detailed implementation settings are provided in Section~\ref{appendix:implementation}.

\subsection{Dataset Preparation for Learning-Based Baselines}
\label{appendix:baseline_dataset_preparation}

To ensure a controlled comparison, all learning-based methods are trained using
the same pool of small C-family Solomon CVRPTW instances used to train
\textsc{Coupled-CG}. Thus, the methods receive the same instance-level training
information, while retaining their original learning procedures. All evaluation
instances are held out from training.

\paragraph{FFCG.}
We adapt \cite{hu2025ffcg} to our CVRPTW column-generation implementation.
The policy interacts with the restricted master problem and pricing procedure
on the training instances, and the resulting transitions are stored in a replay
buffer for DQN training with a target network, following the original learning
framework.

\paragraph{Yuan CG.}
We adapt the PPO-based multiple-column selection method of~\cite{Yuan2024} to the same CVRPTW environment.
The policy is trained on-policy through interaction with column generation.
At each iteration, it observes the generated candidate columns and selects
columns to enter the restricted master problem. For consistency with the main
experimental setting, we use a candidate pool of size $M=10$ and select
$k=5$ columns.

\paragraph{Morabit.}
For \cite{Morabit2021}, we construct the supervised training
data from the same training instances. An expert column-generation procedure
is used to generate column-selection labels and the corresponding bipartite
graph representations, which are then used to train the supervised classifier.

\paragraph{Summary.}
All learning-based baselines therefore use the same CVRPTW training-instance
pool as \textsc{Coupled-CG}, while preserving their respective learning
paradigms: replay-based DQN training for FFCG, on-policy PPO training for
Yuan CG, and expert-labeled supervised learning for Morabit.

\subsection{Implementation Details}
\label{appendix:implementation}

This section reports the architecture, optimization settings, and evaluation
protocol used in our experiments. All methods share the same column-generation
backbone, including the RMP formulation, pricing formulation, exact-pricing fallback, termination criterion, and final integer-master procedure. They differ only in the extension-ordering and/or column-selection policies.
\subsubsection{Offline Training Setup}

Both learning branches are trained using the offline trajectories described in
Section~\ref{appendix:offline_dataset}. Complete episodes are split at the
episode level into $90\%$ training and $10\%$ validation data to avoid leakage
across transitions. Feature-normalization statistics are computed using the
training split only and are fixed during validation and inference. Both the
training and evaluation instances are sampled from the original 100-customer
Solomon CVRPTW instances by constructing smaller subinstances. The training and
test sets are sampled independently, with no overlap between their instances.

We use a batch size of $32$. At each pricing call, at most $M=10$ candidate
columns are retained and at most $k=5$ columns are inserted into the RMP.
When fewer than $k$ valid candidates are available, all valid candidates are
used.

\subsubsection{Shared RMP State Encoder}

The current RMP is represented by the bipartite graph defined in
Section~\ref{appendix:node_feature}, with column nodes, constraint nodes, and
RMP-incidence edges. The shared encoder
$\operatorname{GNN}_{\theta}$ is implemented as a two-layer bipartite
message-passing network with hidden dimension $256$. It outputs a
$128$-dimensional graph representation
\begin{equation}
h_t
=
\operatorname{GNN}_{\theta}(s_t^P),
\end{equation}
and uses dropout rate $0.1$. The same RMP encoder is shared by the
extension-ordering and column-selection branches.

The extension-ordering head maps $h_t$ to five action values corresponding to
the five rules in $\mathcal G$. The prototype-based support estimator described
in Section~\ref{appendix:support_kernel} maintains one prototype for each
rule. We use prototype update rate $\rho=0.05$, Gaussian bandwidth
$\sigma=1.0$, and support-adjustment coefficient $\kappa_P=1.0$.

\subsubsection{Post-Pricing and Relation-Aware Column Encoder}
\label{appendix:column_network}

After pricing, improving routes are sorted by reduced cost, with the
minimum reduced cost route placed first. Each candidate is represented by seven
features: reduced cost, route cost, number of customers, total demand, route
duration, normalized candidate rank, and overlap with the current RMP.

Each candidate feature vector is projected to a $128$-dimensional embedding.
A masked mean of the valid candidate embeddings is concatenated with the shared
RMP representation $h_t$ and transformed into the post-pricing representation
$\widetilde h_t\in\mathbb R^{128}$.

To model interactions among candidate columns, we construct a complete relation
graph over the valid candidates. Each edge $(i,j)$ contains four features:
customer overlap, arc overlap, cosine distance, and Jaccard distance,
\begin{equation}
e_{ij}
=
\left[
e_{ij}^{\mathrm{cust}},
e_{ij}^{\mathrm{arc}},
e_{ij}^{\mathrm{cos}},
e_{ij}^{\mathrm{Jac}}
\right].
\end{equation}
Self-edges and pairs involving padded candidates are masked.

The relation-aware encoder contains two edge-aware attention layers with four
attention heads. For each layer, the relation features contribute to both the
attention score and the transmitted message,
\begin{align}
a_{ij}^{(\ell)}
&=
\frac{
\left(W_Qr_i^{(\ell)}\right)^\top
\left(W_Kr_j^{(\ell)}\right)
}{
\sqrt{d_{\mathrm{head}}}
}
+
w_E^\top e_{ij},\\
m_{ij}^{(\ell)}
&=
W_Vr_j^{(\ell)}
+
W_{EV}e_{ij}.
\end{align}
Masked attention is applied over valid candidates, followed by residual
connections, dropout, and layer normalization.

The resulting candidate representations are passed to the unary and symmetric
pairwise value heads used by the relation-conditioned marginal policy described
in Section~\ref{sec:stage_specific_policies}. Column selection then follows the
marginal procedure defined there; no explicit subset-action catalogue is
constructed.

\subsubsection{Optimization and Training Schedule}

Table~\ref{tab:offline_training_parameters} summarizes the optimization
settings.

\begin{table}[!htbp]
\centering
\caption{Default parameter settings used for offline training.}
\label{tab:offline_training_parameters}
\begin{tabular}{lcc}
\toprule
\textbf{Parameter} & \textbf{Symbol} & \textbf{Value} \\
\midrule
Batch size & $B$ & $32$ \\
Learning rate & -- & $10^{-4}$ \\
Weight decay & -- & $10^{-5}$ \\
Gradient-norm clipping & -- & $5$ \\
Shared embedding dimension & $d_h$ & $128$ \\
GNN hidden dimension & -- & $256$ \\
GNN message-passing layers & -- & $2$ \\
Relation-attention layers & -- & $2$ \\
Relation-attention heads & -- & $4$ \\
Dropout rate & -- & $0.1$ \\
Prototype update rate & $\rho$ & $0.05$ \\
Prototype bandwidth & $\sigma$ & $1.0$ \\
Support coefficient & $\kappa_P$ & $1.0$ \\
Extension conservative coefficient & $\lambda_P$ & $1.0$ \\
Column conservative coefficient & $\lambda_C$ & $1.0$ \\
Polyak update rate & $\tau$ & $0.005$ \\
Optimizer & -- & AdamW \\
Heuristic labeling budget & $\Gamma$ & 200,000\\
\bottomrule
\end{tabular}
\end{table}

Training is performed for $60$ epochs in three stages. During epochs
$1$--$20$, the shared encoder and both decision branches are trained jointly
using alternating updates. During epochs $21$--$40$, the shared encoder and
column-selection branch are frozen and only the extension-ordering branch is
fine-tuned. During epochs $41$--$60$, the shared encoder and
extension-ordering branch are frozen and the post-pricing and column-selection
networks are fine-tuned.

Model selection is performed using held-out validation instances only; test
instances are not used for training or model selection.

\subsubsection{Common Evaluation Protocol}

All comparison methods are evaluated using the same RMP and pricing
implementations, candidate-generation procedure, feasibility and reduced-cost
filters, exact-pricing fallback, termination criterion, candidate-pool size
$M=10$, insertion limit $k=5$, and time limits=$900s$. Thus, the compared methods
differ only in the learned or prescribed decisions being evaluated.

We report both CG time and end-to-end wall-clock time. CG time includes the
RMP and pricing computations,
\begin{equation}
T_{\mathrm{CG}}
=
T_{\mathrm{RMP}}
+
T_{\mathrm{pricing}},
\end{equation}
whereas end-to-end time additionally includes neural inference, process
communication, and orchestration overhead. Completion rate is computed over
the full intended test set.

\subsection{Comparison Methods}
\label{appendix:comparison_methods}

We compare \textsc{Coupled-CG} with four learning-based or heuristic baselines and five ablation configurations. Unless otherwise stated, the candidate-pool size is $M=10$; fixed-cardinality methods insert $k=5$ columns, whereas FFCG and Morabit use variable-size selection. 

\paragraph{Plain CG~\citep{desaulniers2006column}.}
Plain CG uses the fixed minimum reduced cost extension-ordering rule and inserts the $k$ most-negative-reduced-cost columns. No learned component is used.

\paragraph{FFCG~\citep{hu2025ffcg}.}
FFCG uses fixed minimum reduced cost pricing and a DQN-based sequential column selector. Following the original algorithm, we do not impose the fixed cardinality $k$ used by \textsc{Coupled-CG}; candidates are selected one at a time according to marginal Q-values until a learned \textsc{STOP} action is chosen, yielding a variable-size selected set.

\paragraph{Yuan CG~\citep{Yuan2024}.}
Yuan CG uses fixed minimum reduced cost pricing and an on-policy PPO column selector. With $M=10$, the most-negative-reduced-cost column is always included and the policy selects four additional columns from the remaining nine candidates, giving $\binom{9}{4}=126$ subset actions and a fixed $k=5$.

\paragraph{Morabit CG~\citep{Morabit2021}.}
Morabit CG uses fixed minimum reduced cost pricing and a supervised bipartite-graph classifier. Following the original algorithm, we do not restrict the number of inserted columns to $k$. All candidates whose predicted selection probability exceeds $0.5$ are inserted; if none are selected, the most-negative-reduced-cost column is used as a fallback. The selected-set size is therefore variable.

\paragraph{\textsc{Coupled-CG}.}
\textsc{Coupled-CG} jointly uses the learned extension-ordering and relation-aware column-selection policies. The extension policy selects one of five pricing priorities, and the resulting candidate pool is then processed by the column-selection policy. Both policies are trained with the coupled Bellman formulation so that each stage is optimized with respect to the downstream decision.

\paragraph{Pricing Only.}
Pricing Only uses an extension-ordering policy trained independently, without the coupled Bellman updates. The learned pricing policy selects the extension-ordering rule, while column selection follows Plain CG by inserting the $k$ most-negative-reduced-cost columns.

\paragraph{Column Selection Only.}
Column Selection Only uses a column-selection policy trained independently, without coupled Bellman updates from the pricing stage. Pricing is fixed to the minimum reduced cost rule, while the independently trained relation-aware selector chooses the $k$ columns to insert.

\paragraph{Pricing + FFCG.}
This configuration combines the independently trained extension-ordering policy with the FFCG selector. The pricing policy determines the candidate pool, after which FFCG performs variable-size sequential selection.

\paragraph{Pricing + Yuan.}
This configuration combines the independently trained extension-ordering policy with the Yuan PPO selector. The learned pricing policy determines the candidate pool, and Yuan selects a fixed five-column subset.

\paragraph{Pricing + Morabit.}
This configuration combines the independently trained extension-ordering policy with the Morabit classifier. The learned pricing policy determines the candidate pool, after which Morabit performs threshold-based variable-size selection.

The Pricing Only and Column Selection Only policies are trained separately and do not use the cross-stage Bellman coupling of \textsc{Coupled-CG}. The three hybrid configurations retain the same independently trained pricing policy while replacing only the downstream column selector, allowing us to examine how different selectors respond to candidate pools induced by learned pricing.

\subsection{Decoupled Training for Component Comparisons}
\label{appendix:decoupled_training}

To distinguish cross-stage value propagation from independent within-stage learning, the pricing-only and column-selection-only variants, the pricing-only and
column-selection-only variants retain the corresponding state representations,
architectures, rewards, and conservative regularization used in
\textsc{Coupled-CG}, but replace the cross-stage Bellman targets with
within-stage targets. All value-based variants use the undiscounted setting
$\gamma=1$.

\paragraph{Independent pricing model.}
The pricing reward is unchanged:
\begin{equation}
r_t^P
=
-\frac{N_t}{N_{\mathrm{scale}}},
\end{equation}
where $N_t$ is the number of labels explored at iteration $t$ and
$N_{\mathrm{scale}}>0$ is the fixed normalization constant used in
\textsc{Coupled-CG}. The target pricing value at the next pre-pricing state is
\begin{equation}
V_{\bar\psi}^{P}(h_{t+1})
=
\max_{g\in\mathcal G}
\left[
Q_{\bar\psi}^{P}(h_{t+1},g)
-
\kappa_P u_P(h_{t+1},g)
\right].
\end{equation}
Rather than bootstrapping from the column-selection value at the same
iteration, the independent pricing model bootstraps directly from the next
pricing state:
\begin{equation}
y_{t,\mathrm{ind}}^{P}
=
r_t^P
+
(1-\mathrm{done}_t)
V_{\bar\psi}^{P}(h_{t+1}).
\label{eq:ind_pricing_target}
\end{equation}
Thus, pricing value propagates only across consecutive pricing states,
\[
Q_t^P
\longleftarrow
Q_{t+1}^P,
\]
and does not directly bootstrap from $Q_t^C$.

\paragraph{Independent column-selection model.}
The column-selection reward is also unchanged:
\begin{equation}
r_t^C
=
\beta
\frac{Z_t-Z_{t+1}}
{\max\{1,|Z_1|\}}.
\end{equation}
For the next post-pricing state
$\widetilde h_{t+1}
=
\operatorname{Enc}_{\theta}^{C}(h_{t+1},\mathcal C_{t+1})$,
the target network greedily constructs
\begin{equation}
\widehat O_{t+1}
=
\operatorname{GreedySelect}_{\bar\phi}
\left(
Q_{\bar\phi}^{C},
\widetilde h_{t+1}
\right),
\end{equation}
and evaluates
\begin{equation}
V_{\bar\phi}^{C}(\widetilde h_{t+1})
=
Q_{\bar\phi}^{C}
\left(
\widetilde h_{t+1},
\widehat O_{t+1}
\right).
\end{equation}
The independent column-selection target is therefore
\begin{equation}
y_{t,\mathrm{ind}}^{C}
=
r_t^C
+
(1-\mathrm{done}_t)
V_{\bar\phi}^{C}(\widetilde h_{t+1}).
\label{eq:ind_column_target}
\end{equation}
Hence, column-selection value propagates only across consecutive
column-selection states,
\[
Q_t^C
\longleftarrow
Q_{t+1}^C,
\]
rather than bootstrapping from the next pricing value $Q_{t+1}^P$.

\paragraph{Difference from coupled training.}
The independent variants therefore use the within-stage propagation
\begin{equation}
Q_t^P
\longleftarrow
Q_{t+1}^P,
\qquad
Q_t^C
\longleftarrow
Q_{t+1}^C,
\end{equation}
whereas \textsc{Coupled-CG} uses the cross-stage Bellman structure
\begin{equation}
Q_t^P
\longleftarrow
Q_t^C
\longleftarrow
Q_{t+1}^P.
\end{equation}
In the independent variants, the other decision stage remains part of the
logged CG transition dynamics, but its learned value is not used in the
Bellman backup. Since the independent and coupled variants use the same
undiscounted setting $\gamma=1$, rewards, and corresponding model components,
their distinction is the presence or absence of cross-stage value propagation.

\subsection{Computational Environment}
\label{appendix:computational_environment}

All experiments are conducted on a Linux server running Ubuntu 22.04-compatible
Linux (kernel 5.15, x86-64). The machine is equipped with two Intel Xeon Gold
5220 processors, providing 36 physical cores and 72 hardware threads in total,
with a base frequency of 2.20\,GHz and a maximum frequency of 3.90\,GHz. The
system contains 188\,GiB of RAM and one NVIDIA Quadro RTX 6000 GPU with
24\,GB of VRAM.

The learning components are implemented in Python using PyTorch with CUDA 12.1.
The main Python dependencies include \texttt{numpy}, \texttt{torch},
\texttt{gurobipy}, \texttt{pybind11}, \texttt{scikit-build-core},
\texttt{ninja}, and \texttt{PyYAML}. The native components are compiled with
GCC/G++ 11.4, CMake 3.22.1, and Ninja 1.10.1.

The RMP and other mathematical optimization components are solved using
Gurobi through \texttt{gurobipy}. All experiments use the same computational
environment and solver configuration across comparison methods to ensure a
controlled runtime comparison.

\subsection{Overall Performance and Generalization}
\paragraph{Overall performance.}
Across all four evaluation groups, Table~\ref{tab:small_large_results} shows that \textsc{Coupled-CG} consistently lowers computational cost. Relative to Plain CG, CG time decreases by 62.5\%, 43.1\%, 76.6\%, and 39.4\% on C-Small, C-Large, R/RC-Small, and R/RC-Large, respectively, with corresponding reductions in explored labels of 67.7\%, 60.2\%, 78.2\%, and 44.2\%. On the instances solved by all five methods, \textsc{Coupled-CG} achieves the lowest mean CG and end-to-end runtime in every group. CG Solved (\%) is identical across methods in the two small-instance groups and highest for \textsc{Coupled-CG} on C-Large, while on R/RC-Large it is slightly lower than for Yuan and Morabit; the reported runtime gains therefore characterize the jointly solved subset rather than uniform dominance over the full test set.

\paragraph{Generalization.}
Because all learning-based methods are trained only on small C-family instances and every evaluation instance is held out, the four groups probe different levels of transfer. C-Small evaluates unseen instances within the training family and size range, C-Large tests generalization to larger problem sizes, R/RC-Small tests transfer to different instance structures, and R/RC-Large combines both shifts. \textsc{Coupled-CG} continues to achieve substantial reductions in CG time and labeling effort across all of these settings, including the larger R/RC cases, showing that the learned guidance transfers beyond the training distribution.
\rev{Beyond CVRPTW (Section~\ref{appendix:generalization_beyond_cvrptw}), the CVRPTW-trained policy, applied without retraining to VRPTW without capacity, explores 64\% fewer labels and needs 69\% less CG time than Plain CG on large instances; trained for airline crew pairing, \textsc{Coupled-CG} reduces labels by 41--55\% and CG time by 35--45\%.}

\begin{table}[!htbp]
\centering
\caption{Performance comparison on small and large instances. Runtime, labeling effort, and CG iterations are averaged over instances solved by all five methods within the 900-second time limit, while \textbf{CG Solved (\%)} reports the percentage of all instances in each group solved by each method within the time limit.}
\label{tab:small_large_results}

\small
\setlength{\tabcolsep}{3.5pt}

\resizebox{\textwidth}{!}{%
\begin{tabular}{clrrrrr}
\toprule
\textbf{Group} & \textbf{Method}
& \textbf{CG Time (s)}
& \textbf{End-to-End (s)}
& \textbf{Labels Explored}
& \textbf{CG Iterations}
& \textbf{CG Solved (\%)} \\
\midrule

\multirow{5}{*}{\makecell{ C-Small\\(6--18 customers)\\Solved by all: 281}}
& \textsc{Coupled-CG} & \textbf{2.12} & \textbf{2.15}
& \textbf{22,662} & \textbf{18.33} & 99.65 \\

& Plain CG~\citep{desaulniers2006column} & 5.66 & 5.69
& 70,091 & 19.17 & 99.65 \\

& Yuan CG~\citep{Yuan2024} & 6.08 & 6.25
& 72,759 & 19.44 & 99.65 \\

& Morabit~\citep{Morabit2021} & 5.70 & 5.70
& 118,104 & 42.79 & 99.65 \\

& FFCG~\citep{hu2025ffcg} & 3.87 & 3.87
& 66,925 & 20.32 & 99.65 \\

\midrule

\multirow{5}{*}{\makecell{C-Large\\(20--60 customers)\\Solved by all: 256}}
& \textsc{Coupled-CG} & \textbf{29.54} & \textbf{29.56}
& \textbf{271,943} & \textbf{86.03} & \textbf{84.47} \\

& Plain CG~\citep{desaulniers2006column} & 51.92 & 51.97
& 683,037 & 102.79 & 83.23 \\

& Yuan CG~\citep{Yuan2024} & 50.20 & 50.40
& 678,558 & 102.38 & 83.54 \\

& Morabit~\citep{Morabit2021} & 44.55 & 44.55
& 947,024 & 91.08 & 82.30 \\

& FFCG~\citep{hu2025ffcg} & 44.50 & 44.51
& 863,216 & 93.14 & 83.23 \\

\midrule

\multirow{5}{*}{\makecell{R/RC-Small\\ (6--18 customers)\\Solved by all: 186}}
& \textsc{Coupled-CG} & \textbf{2.90} & \textbf{3.29}
& \textbf{27,399} & 15.90 & 85.32 \\

& Plain CG~\citep{desaulniers2006column} & 12.37 & 12.39
& 125,670 & 16.52 & 85.32 \\

& Yuan CG~\citep{Yuan2024} & 12.82 & 12.95
& 120,832 & 16.19 & 85.32 \\

& Morabit~\citep{Morabit2021} & 6.68 & 6.68
& 93,110 & \textbf{12.19} & 85.32 \\

& FFCG~\citep{hu2025ffcg} & 6.85 & 6.86
& 101,611 & 17.04 & 85.32 \\

\midrule

\multirow{5}{*}{\makecell{R/RC-Large \\(20--60 customers)\\Solved by all: 200}}
& \textsc{Coupled-CG} & \textbf{21.66} & \textbf{21.70}
& \textbf{191,331} & \textbf{61.16} & 74.46 \\

& Plain CG~\citep{desaulniers2006column} & 35.74 & 35.86
& 342,953 & 67.07 & 72.30 \\

& Yuan CG~\citep{Yuan2024} & 29.62 & 29.93
& 343,698 & 67.21 & 75.90 \\

& Morabit~\citep{Morabit2021} & 26.77 & 27.78
& 269,481 & 65.43 & \textbf{76.26} \\

& FFCG~\citep{hu2025ffcg} & 29.43 & 29.44
& 362,566 & 61.61 & 74.82 \\

\bottomrule
\end{tabular}%
}
\end{table}

\subsection{Component and Hybrid Comparisons}
\begin{figure}[!htbp]
    \centering
    \includegraphics[width=\linewidth]{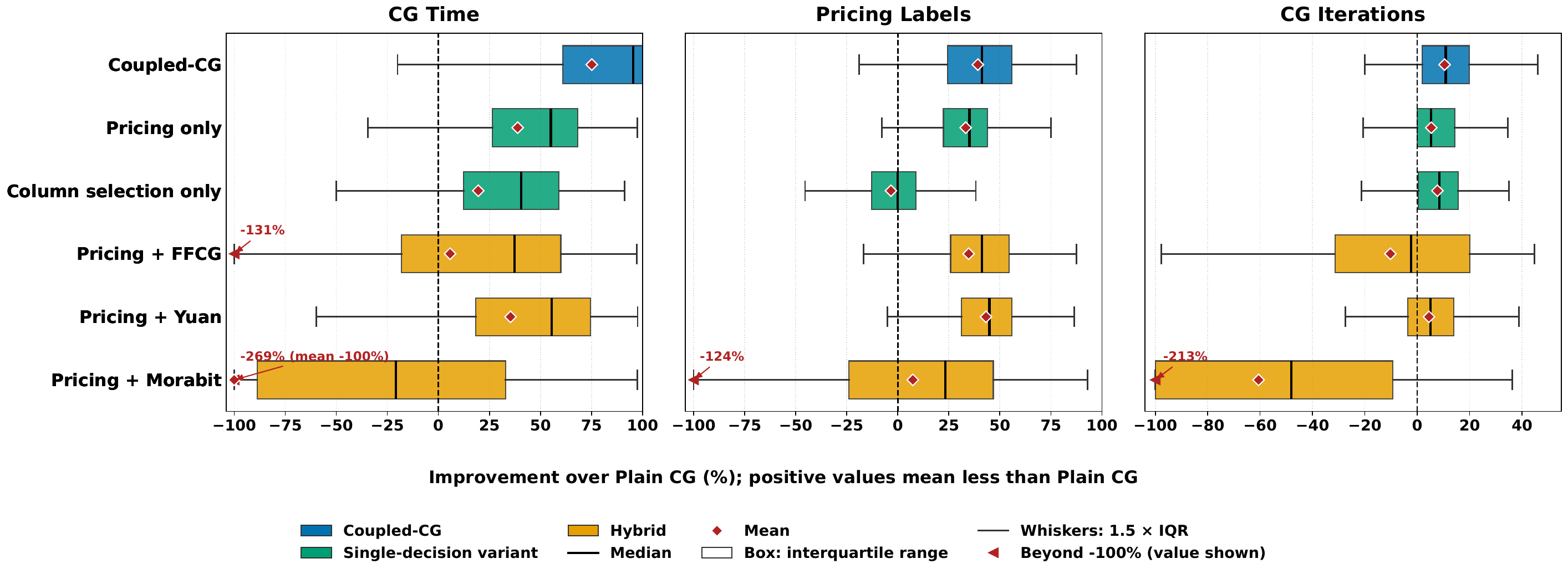}
    \caption{Ablation study of pricing guidance and column selection. Positive
    values indicate improvement over Plain CG.}
    \label{fig:ablation}
\end{figure}

Figure~\ref{fig:ablation} compares the contributions of learned pricing guidance and column selection and examines the additional benefit of their coupled training. The pricing-only and column-selection-only variants retain the corresponding architectures and rewards of \textsc{Coupled-CG}, but are trained independently using within-stage TD targets (Details in Section~\ref{appendix:decoupled_training}). All three value-based variants use the undiscounted setting \(\gamma=1\), eliminating discounting as a source of variation across these comparisons. Both independent components improve over Plain CG, while \textsc{Coupled-CG} yields the most consistent gains in CG time, pricing effort, and CG iterations. Combining the independently trained pricing policy with FFCG, Yuan, or Morabit produces substantially larger variance and, in some cases, severe degradation, suggesting sensitivity to the candidate-column distribution induced by learned pricing.
\rev{We also ablate the two refinements within \textsc{Coupled-CG} (Table~\ref{tab:refinement_ablation}), replacing support-aware ordering with a standard CQL argmax and relation-conditioned selection with per-column top-$k$ scoring. Only the full model solves all 300 instances; removing either refinement raises mean labels by 73\% and 432\%, respectively.}

\begin{table}[!htbp]
\centering
\caption{\rev{Ablation of the two refinements (300 small test instances). Means are over solved instances.}}
\label{tab:refinement_ablation}
\small
\begin{tabular}{cccccc}
\toprule
\textbf{Support-aware} & \textbf{Relation-cond.} & \textbf{Exact term.} & \textbf{Iterations} & \textbf{Labels} & \textbf{Time (s)} \\
\midrule
\checkmark & \checkmark & 300/300 & 14.8 & 21{,}895 & 2.34 \\
\checkmark & -- & 274/300 & 16.95 & 37{,}826 & 1.20 \\
-- & \checkmark & 266/300 & 17.43 & 116{,}535 & 5.42 \\
-- & -- & 260/300 & 17.05 & 119{,}189 & 6.44 \\
\bottomrule
\end{tabular}
\end{table}

\subsection{Sensitivity to Candidate-Pool Size and Selection Cardinality}

We further evaluate the sensitivity of \textsc{Coupled-CG} to the candidate-pool size $M$ and selection cardinality $k$ on C-family instances with 40--60 customers. Table~\ref{tab:c_40_60_mk} evaluates whether the advantage of \textsc{Coupled-CG} depends on the specific candidate-pool and selection sizes used in the main experiments. Across all six $M/k$ settings, \textsc{Coupled-CG} consistently reduces CG time, end-to-end time, labeling effort, and the number of CG iterations relative to Plain CG. Although the magnitude of the improvement varies with $M$ and $k$, the gains persist from the smallest setting $(10,5)$ to the largest $(60,30)$. This indicates that the benefit of the proposed framework is not tied to a particular candidate-set configuration.

\begin{table}[!htbp]
\centering
\caption{Performance comparison under different $M/k$ settings on C instances with 40-60 customers. All metrics are averaged over the instances solved by both \textsc{Coupled-CG} and Plain CG within the time limit.}
\label{tab:c_40_60_mk}
\resizebox{\textwidth}{!}{
\begin{tabular}{llrrrrr}
\toprule
\textbf{$M/k$}
& \textbf{Method}
& \textbf{Solved by Both}
& \textbf{CG Time (s)}
& \textbf{End-to-End (s)}
& \textbf{Labels Explored}
& \textbf{CG Iterations} \\
\midrule

\multirow{2}{*}{$10/5$}
& \textsc{Coupled-CG}
& 22
& \textbf{43.61}
& \textbf{45.61}
& \textbf{419,024}
& \textbf{136.59} \\

& Plain CG
& 22
& 81.92
& 82.07
& 885,506
& 145.73 \\

\midrule

\multirow{2}{*}{$20/10$}
& \textsc{Coupled-CG}
& 25
& \textbf{53.43}
& \textbf{55.92}
& \textbf{403,916}
& \textbf{96.20} \\

& Plain CG
& 25
& 132.37
& 132.78
& 996,610
& 110.44 \\

\midrule

\multirow{2}{*}{$30/15$}
& \textsc{Coupled-CG}
& 25
& \textbf{106.11}
& \textbf{108.79}
& \textbf{433,809}
& \textbf{80.04} \\

& Plain CG
& 25
& 120.21
& 120.48
& 767,626
& 87.68 \\

\midrule

\multirow{2}{*}{$40/20$}
& \textsc{Coupled-CG}
& 26
& \textbf{96.67}
& \textbf{99.26}
& \textbf{390,730}
& \textbf{69.04} \\

& Plain CG
& 26
& 144.56
& 144.82
& 800,774
& 78.58 \\

\midrule

\multirow{2}{*}{$50/25$}
& \textsc{Coupled-CG}
& 26
& \textbf{80.06}
& \textbf{83.12}
& \textbf{352,220}
& \textbf{65.04} \\

& Plain CG
& 26
& 106.03
& 106.31
& 727,022
& 74.35 \\

\midrule

\multirow{2}{*}{$60/30$}
& \textsc{Coupled-CG}
& 24
& \textbf{53.31}
& \textbf{56.15}
& \textbf{272,305}
& \textbf{57.63} \\

& Plain CG
& 24
& 87.67
& 87.99
& 639,016
& 68.96 \\

\bottomrule
\end{tabular}
}
\end{table}

\subsection{Which Instance Structures Benefit Most from COUPLED-CG?}
\begin{figure}[!htbp]
    \centering
    \includegraphics[width=\textwidth]{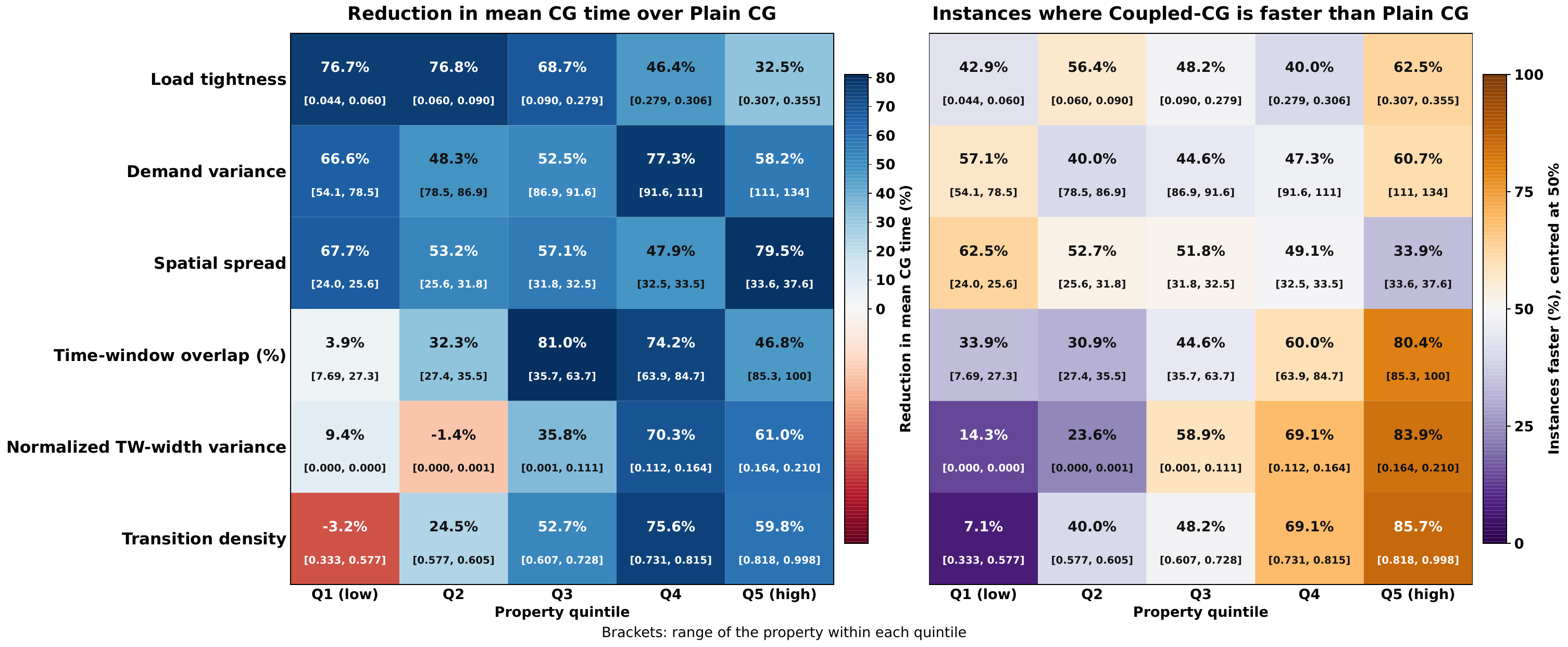}
    \caption{Structural analysis of \textsc{Coupled-CG} across CVRPTW instance properties. The left panel reports the reduction in mean CG time relative to Plain CG,
$100\left(1-\bar T^{\mathrm{Coupled}}/\bar T^{\mathrm{Plain}}\right)$,
within each property quintile, while the right panel reports the percentage of instances for which \textsc{Coupled-CG} is faster than Plain CG.}
    \label{fig:structural_pca}
\end{figure}

We investigate which instance structures benefit most from \textsc{Coupled-CG}. Each instance is characterized by six properties: load tightness, demand variance, spatial spread, time-window overlap, normalized time-window-width variance, and transition density; detailed definitions are provided in Section~\ref{appendix:structural_features}. For each property, we partition instances into quintiles, and Figure~\ref{fig:structural_pca} reports both the average CG-time reduction over Plain CG and the percentage of instances where \textsc{Coupled-CG} is faster.

The clearest trends arise from temporal complexity and transition density. As transition density increases from Q1 to Q5, the average CG-time reduction rises from $-3.2\%$ to $59.8\%$, while the fraction of instances where \textsc{Coupled CG} is faster increases from $7.1\%$ to $85.7\%$. Similar improvements in consistency are observed for time-window overlap and normalized time-window-width variance, whereas load tightness, demand variance, and spatial spread show less monotonic effects. These results indicate that \textsc{Coupled-CG} is most promising for routing problems with dense feasible transitions and heterogeneous temporal constraints, such as ride-pooling, pickup-and-delivery, and time-window-constrained last-mile delivery.

Additional sensitivity analyses with the reward coefficient $\beta$ are provided in
Section~\ref{appendix:additional_analysis}.

\subsection{Structural Instance Features}
\label{appendix:structural_features}

We characterize each CVRPTW instance using six structural properties capturing
capacity, demand, spatial, temporal, and transition characteristics. Let
$\mathcal N$ denote the customer set, $n=|\mathcal N|$, $q_i$ the demand of
customer $i$, $Q$ the vehicle capacity, $K$ the number of vehicles,
$(x_i,y_i)$ the customer coordinates, $[a_i,b_i]$ the time window, and $s_i$
the service time.

\paragraph{Load tightness.}
Load tightness measures total demand relative to the available fleet capacity:
\begin{equation}
\mathrm{LT}
=
\frac{\sum_{i\in\mathcal N} q_i}{KQ}.
\end{equation}
A larger value indicates tighter aggregate vehicle capacity.

\paragraph{Demand variance.}
Demand variance measures the heterogeneity of customer demands:
\begin{equation}
\mathrm{DV}
=
\frac{1}{n}
\sum_{i\in\mathcal N}(q_i-\bar q)^2,
\qquad
\bar q=\frac{1}{n}\sum_{i\in\mathcal N}q_i.
\end{equation}
Larger values indicate greater variation in customer load requirements.

\paragraph{Spatial spread.}
Spatial spread measures the dispersion of customer locations along both spatial
dimensions:
\begin{equation}
\mathrm{SS}
=
\sqrt{
\operatorname{Var}(x)
+
\operatorname{Var}(y)
}.
\end{equation}
A larger value indicates a more geographically dispersed customer set.

\paragraph{Time-window overlap.}
For each unordered customer pair $(i,j)$, we identify whether their time
windows have a nonempty intersection:
\begin{equation}
\mathbb{I}_{ij}^{\mathrm{TW}}
=
\mathbb{I}
[
\max(a_i,a_j)
\leq
\min(b_i,b_j)
].
\end{equation}
The instance-level overlap percentage is
\begin{equation}
\mathrm{TWO}
=
\frac{100}{\binom{n}{2}}
\sum_{i<j}
\mathbb{I}_{ij}^{\mathrm{TW}}.
\end{equation}
A larger value indicates stronger overlap among customer service windows.

\paragraph{Normalized time-window-width variance.}
Let $H=b_0-a_0$ denote the depot time horizon and
$w_i=b_i-a_i$ the width of customer $i$'s time window. We first normalize
each width as
\begin{equation}
\widetilde w_i
=
\frac{w_i}{H},
\end{equation}
and define
\begin{equation}
\mathrm{TWV}
=
\frac{1}{n}
\sum_{i\in\mathcal N}
(\widetilde w_i-\bar{\widetilde w})^2.
\end{equation}
Larger values indicate greater heterogeneity in relative time-window widths.

\paragraph{Transition density.}
We construct a directed feasible-transition graph over customers. Let $d_{0i}$
denote the depot-to-customer distance and $d_{ij}$ the distance from customer
$i$ to $j$. The earliest feasible service time at customer $i$ is
\begin{equation}
e_i
=
\max\{a_i,a_0+d_{0i}\}.
\end{equation}
A directed transition $i\rightarrow j$ is considered feasible if
\begin{equation}
q_i+q_j\leq Q,
\qquad
e_i+s_i+d_{ij}\leq b_j,
\end{equation}
and customer $j$ can subsequently return to the depot within the depot horizon:
\begin{equation}
\max\{a_j,e_i+s_i+d_{ij}\}
+s_j+d_{j0}
\leq b_0.
\end{equation}
Let $A_{ij}=1$ when these conditions hold and $A_{ij}=0$ otherwise. The
transition density is
\begin{equation}
\mathrm{TD}
=
\frac{1}{n(n-1)}
\sum_{i\neq j} A_{ij}.
\end{equation}
A larger value indicates a denser feasible-transition graph and more possible
customer-to-customer continuations.

\subsection{Sensitivity to the Reward Coefficient}\label{appendix:additional_analysis}

Table~\ref{tab:beta_sensitivity} examines the effect of the reward coefficient $\beta$. The main setting $\beta=1.0$ was fixed before this sweep and was not selected retrospectively based on the best observed performance. Across $\beta\in[0.1,1.0]$, CG time remains consistently lower than the $\beta=0$ setting, with reductions ranging from approximately $4\%$ to $9\%$. The remaining metrics vary only moderately across the tested values, suggesting that the method is not highly sensitive to the precise choice of $\beta$ within this range. In particular, $\beta=1.0$ achieves performance comparable to the best settings while preserving the reward formulation used throughout the main experiments.

\begin{table}[!htbp]
\centering
\caption{Sensitivity analysis of the reward coefficient $\beta$.
The setting $\beta=1.0$, used in the main experiments, is marked in bold.
The comparison is reported for sensitivity analysis only; $\beta$ was not
selected retrospectively based on the best performance in this experiment.}
\label{tab:beta_sensitivity}
\resizebox{\textwidth}{!}{
\begin{tabular}{rrrrrrr}
\toprule
\textbf{$\beta$}
& \textbf{Total CG Time (s)}
& \textbf{Mean CG Time (s)}
& \textbf{vs.\ $\beta=0$}
& \textbf{Wall Time (s)}
& \textbf{Labels}
& \textbf{Iterations} \\
\midrule

0.0
& 189.14
& 3.940
& -- 
& 296.22
& 99,687
& 493 \\

0.1
& 179.44
& 3.738
& $-5.1\%$
& 284.87
& 101,005
& 499 \\

0.2
& 181.66
& 3.785
& $-4.0\%$
& 286.90
& 101,575
& 499 \\

0.3
& 173.75
& 3.620
& $-8.1\%$
& 277.58
& 102,060
& 496 \\

0.4
& 175.16
& 3.649
& $-7.4\%$
& 279.07
& 101,360
& 497 \\

0.5
& 173.52
& 3.615
& $-8.3\%$
& 277.65
& \textbf{99,606}
& \textbf{493} \\

0.6
& 175.97
& 3.666
& $-7.0\%$
& 280.44
& 101,957
& 499 \\

0.7
& 174.42
& 3.634
& $-7.8\%$
& 278.62
& 99,735
& \textbf{493} \\

0.8
& \textbf{172.31}
& \textbf{3.590}
& \textbf{$-8.9\%$}
& 277.63
& 101,549
& 496 \\

0.9
& 178.39
& 3.716
& $-5.7\%$
& 282.94
& 101,282
& 496 \\

\textbf{1.0}
& 173.82
& 3.621
& $-8.1\%$
& 278.31
& 102,008
& 496 \\

\bottomrule
\end{tabular}
}
\end{table}

\subsection{\rev{Generalization Beyond CVRPTW}}
\label{appendix:generalization_beyond_cvrptw}

\rev{\textsc{Coupled-CG} relies on two properties of the CG algorithm, not on CVRPTW itself: pricing is a resource-constrained shortest path solved by labeling, and several candidate columns are available for insertion at each iteration. We test it on two further problems with these properties. In both, the reference is Plain CG, which uses the minimum-reduced-cost extension order and inserts the $k$ most-negative-reduced-cost columns, with $M=10$ and $k=5$. Exact pricing certifies termination, so both methods reach the same LP optimum on every instance they solve.}

\subsubsection{\rev{VRPTW Without Capacity}}

\rev{We remove the capacity constraint from CVRPTW and apply the policy trained on CVRPTW without retraining. The test set has 60 small and 60 large instances, 10 per Solomon family, with a 900\,s time limit. Table~\ref{tab:vrptw_no_capacity} reports the ratio of the \textsc{Coupled-CG} mean to the Plain CG mean over the instances that both methods solve within the limit. On the large set, the transferred policy explores 64\% fewer labels and needs 69\% less CG time. On the small set, it explores 28\% fewer labels, but CG time is 7\% higher, so the label reduction does not yield a runtime gain at this size.}

\begin{table}[!htbp]
\centering
\caption{\rev{VRPTW without capacity: ratio of means, \textsc{Coupled-CG} to Plain CG, over commonly solved instances. The policy is trained on CVRPTW and applied without retraining.}}
\label{tab:vrptw_no_capacity}
\begin{tabular}{lcccc}
\toprule
\textbf{Set} & \textbf{Instances} & \textbf{Iterations} & \textbf{Labels} & \textbf{CG time} \\
\midrule
Small & 58/60 & $0.91\times$ & $0.72\times$ & $1.07\times$ \\
Large & 36/60 & $0.95\times$ & $0.36\times$ & $0.31\times$ \\
\bottomrule
\end{tabular}
\end{table}

\subsubsection{\rev{Airline Crew Pairing}}

\rev{The airline crew pairing problem covers every flight leg exactly once with pairings. A pairing is a sequence of duties, separated by overnight rests, that starts and ends at the same crew base. Pricing is a resource-constrained shortest path solved by labeling. We use the instances of \citet{kasirzadeh2017airline}: seven monthly airline schedules with 1{,}013 to 7{,}766 legs, which we cut into day windows. Schedules 1, 4, 6, and 7 are held out for testing, and the policy is trained on windows from the remaining schedules. The test set has 52 small windows (132--261 legs) and 51 large windows (323--1{,}013 legs).}

\rev{The two-stage MDP carries over with problem-specific columns, resources, and ordering rules (Table~\ref{tab:crew_pairing_mdp}). Table~\ref{tab:crew_pairing_results} reports the geometric mean of the \textsc{Coupled-CG} to Plain CG ratio over the test windows, all of which both methods solve. CG time counts pricing and RMP solution time and excludes model inference. \textsc{Coupled-CG} explores 41\% and 55\% fewer labels on the small and large sets, and reduces CG time by 35\% and 45\%.}

\begin{table}[!htbp]
\centering
\caption{\rev{The two-stage MDP for CVRPTW and for airline crew pairing.}}
\label{tab:crew_pairing_mdp}
\small
\begin{tabular}{lll}
\toprule
\textbf{MDP element} & \textbf{CVRPTW} & \textbf{Airline crew pairing} \\
\midrule
Column & Vehicle route & Pairing \\
RMP rows & Customer covering, fleet size & Leg covering \\
Label resources & Reduced cost, load, time & \makecell[l]{Reduced cost, pairing start, duties used,\\ duty start, duty flying time, duty legs} \\
Pricing action $g_t$ & \makecell[l]{Order labels by reduced cost,\\ time, or load} & \makecell[l]{Order labels by reduced cost, legs covered,\\ time away from base, or unproductive time} \\
Selection action $O_t$ & Up to $k$ routes & Up to $k=5$ of $M=10$ pairings \\
\bottomrule
\end{tabular}
\end{table}

\begin{table}[!htbp]
\centering
\caption{\rev{Airline crew pairing: geometric mean of the \textsc{Coupled-CG} to Plain CG ratio on held-out test windows.}}
\label{tab:crew_pairing_results}
\begin{tabular}{lcccc}
\toprule
\textbf{Set} & \textbf{Instances} & \textbf{Iterations} & \textbf{Labels} & \textbf{CG time} \\
\midrule
Small (132--261 legs) & 52/52 & $0.84\times$ & $0.59\times$ & $0.65\times$ \\
Large (323--1{,}013 legs) & 51/51 & $0.82\times$ & $0.45\times$ & $0.55\times$ \\
\bottomrule
\end{tabular}
\end{table}

\subsection{\rev{Training Time}}
\label{appendix:training_time}

\rev{Table~\ref{tab:training_time} reports the training time of the learning-based methods on the same training data, over three runs per method. Among the RL methods, \textsc{Coupled-CG} trains fastest: its mean time is $2.6\times$ lower than that of FFCG and $13.0\times$ lower than that of Yuan CG, although it learns two decisions instead of one. The reason is that it learns from fixed trajectories and never runs CG during training, whereas FFCG and Yuan CG interact with the CG solver. Morabit CG trains faster because it is supervised: it fits a classifier to fixed expert labels, while RL bootstraps long-horizon values and needs repeated target updates. Its time excludes the 0.30\,h needed to generate the expert labels.}

\begin{table}[!htbp]
\centering
\caption{\rev{Training time of the learning-based methods on the same training data (three runs per method).}}
\label{tab:training_time}
\begin{tabular}{lllcc}
\toprule
\textbf{Method} & \textbf{Learning} & \textbf{Decisions learned} & \textbf{Mean (h)} & \textbf{Median (h)} \\
\midrule
Morabit CG & Supervised classifier & Column selection & 0.34 & 0.30 \\
\textsc{Coupled-CG} & Offline RL & Pricing and column selection & 1.44 & 0.91 \\
FFCG & Online RL & Column selection & 3.80 & 3.82 \\
Yuan CG & Online RL & Column selection & 18.82 & 18.61 \\
\bottomrule
\end{tabular}
\end{table}

\subsection{\rev{Comparison with Learning-Based Pricing Methods}}
\label{appendix:learned_pricing_comparison}

\rev{The baselines in Table~\ref{tab:small_large_results} keep pricing exact and differ only in column selection. Learning-based pricing methods instead prune or replace the pricing network, so on their own they do not certify the LP optimum, and their effort is not directly comparable. We make them comparable by adding an exact pricing fallback, so that every method terminates at the same LP optimum. We implement three such methods in the same CG code with $M=10$, $k=5$, and a 900\,s time limit: arc selection~\citep{morabit2023}, RLHH~\citep{Xu2023RLHH}, and POMO-CG~\citep{Abouelrous2025RLPricing}.}

\rev{Table~\ref{tab:learned_pricing_comparison} reports the results on 300 small and 300 large test instances. Means and medians are over the instances that all exact methods solve (299 small, 105 large). On small instances, every learned method reaches exact termination on all 300 instances; arc selection and RLHH have lower mean time than \textsc{Coupled-CG} (0.85\,s and 1.00\,s against 2.34\,s), and \textsc{Coupled-CG} explores the fewest labels. On large instances, the three learning-based pricing methods solve fewer instances than Plain CG (106--126 against 185) and explore 2.5 to 4.9 times as many labels. \textsc{Coupled-CG} solves 239 of 300 large instances, with fewer labels and lower mean time than Plain CG.}

\begin{table}[!htbp]
\centering
\caption{\rev{Comparison with learning-based pricing methods under an exact pricing fallback. Means and medians are over the instances solved by all exact methods (299 small, 105 large).}}
\label{tab:learned_pricing_comparison}
\resizebox{\textwidth}{!}{
\begin{tabular}{llccrrrr}
\toprule
\textbf{Set} & \textbf{Method} & \textbf{Exact termination} & \textbf{Time limit} & \textbf{Mean iterations} & \textbf{Mean labels} & \textbf{Mean time (s)} & \textbf{Median time (s)} \\
\midrule
Small & Plain CG & 299/300 & 1 & 16.8 & 103{,}323 & 12.84 & 0.04 \\
Small & Arc selection & 300/300 & 0 & 16.6 & 22{,}742 & 0.85 & 0.07 \\
Small & RLHH & 300/300 & 0 & 26.7 & 24{,}484 & 1.00 & 0.07 \\
Small & POMO-CG (hybrid) & 300/300 & 0 & 15.6 & 48{,}990 & 7.51 & 0.10 \\
Small & \textsc{Coupled-CG} & 300/300 & 0 & 14.8 & 21{,}895 & 2.34 & 0.06 \\
\midrule
Large & Plain CG & 185/300 & 115 & 50.0 & 119{,}708 & 12.02 & 0.96 \\
Large & Arc selection & 126/300 & 174 & 49.5 & 393{,}524 & 35.38 & 5.35 \\
Large & RLHH & 126/300 & 174 & 70.1 & 300{,}700 & 18.25 & 3.65 \\
Large & POMO-CG (hybrid) & 106/300 & 194 & 49.6 & 585{,}669 & 118.77 & 16.01 \\
Large & \textsc{Coupled-CG} & 239/300 & 61 & 45.2 & 87{,}814 & 8.76 & 1.01 \\
\bottomrule
\end{tabular}
}
\end{table}

\subsection{\rev{Comparison and Combination with Adaptive Stabilization}}
\label{appendix:ascg_comparison}

\rev{ASCG-ML~\citep{Shen2024AdaptiveCG} accelerates CG by dual stabilization. It acts on a different part of a CG iteration than \textsc{Coupled-CG} (Table~\ref{tab:ascg_method}). ASCG-ML changes the dual vector that is sent to pricing. \textsc{Coupled-CG} takes the dual vector it is given and decides how the pricing search is ordered and which returned columns are added. Neither method modifies the component the other controls, so both can run in the same CG loop.}

\begin{table}[!htbp]
\centering
\caption{\rev{ASCG-ML and \textsc{Coupled-CG} act on different parts of a CG iteration.}}
\label{tab:ascg_method}
\small
\begin{tabular}{p{0.16\textwidth}p{0.38\textwidth}p{0.38\textwidth}}
\toprule
 & \textbf{ASCG-ML} & \textbf{\textsc{Coupled-CG}} \\
\midrule
What is learned & A prediction of the optimal dual solution, made once per instance & Two policies applied at every iteration: label-extension ordering and column selection \\
Where it acts & Dual side of the RMP: the predicted duals are the stabilization center, and a penalty keeps the RMP duals close to it & Pricing search and primal side of the RMP: how labeling explores, and which generated columns enter the RMP \\
How it adapts & The penalty weight is updated from the current reduced cost and removed when no improving column is found & Both policies observe the current RMP state (columns and duals) at each iteration \\
Learning & Supervised regression on optimal duals & Offline RL on logged CG trajectories \\
Source of speed-up & Fewer CG iterations, by damping dual oscillation & Less labeling effort per iteration and more useful columns \\
Exactness & Final exact pricing under the unstabilized duals & Final exact pricing \\
\bottomrule
\end{tabular}
\end{table}

\paragraph{\rev{Combined algorithm.}}
\rev{Let $\hat{\pi}$ be the predicted duals and $\rho_t$ the penalty weight. At each iteration $t$, the combined algorithm (\textsc{Coupled-CG} + ASCG-ML) performs five steps. (1) \emph{Stabilized RMP}: solve the RMP with penalty $\rho_t$ around $\hat{\pi}$, and obtain the stabilized duals $\tilde{\pi}_t$ and the objective $Z_t$. (2) \emph{Pricing stage}: encode the RMP state with $\tilde{\pi}_t$, select the extension-ordering rule $g_t$ with $Q^P$, and run labeling under $\tilde{\pi}_t$ to obtain the candidate set $\mathcal C_t$. (3) \emph{Column-selection stage}: select $O_t\subseteq\mathcal C_t$ with $|O_t|\leq k$ using $Q^C$, and add $O_t$ to the RMP. (4) \emph{Penalty update}: update $\rho_{t+1}$ with the ASCG rule from the best reduced cost found in step (2). (5) \emph{Termination}: if $\mathcal C_t=\emptyset$, set the penalty to zero, re-solve the RMP, and run exact pricing under the unstabilized duals; stop if no column with negative reduced cost exists, and otherwise add the columns and continue. Step (5) is the termination test used by each method alone, so the combined algorithm returns the same LP optimum.}

\paragraph{\rev{Results.}}
\rev{We implemented ASCG-ML in our CVRPTW CG code and trained it on our training instances. Table~\ref{tab:ascg_results} reports means over the instances on which all methods finished (289 small, 207 large); all methods reach the same LP optimum. \textsc{Coupled-CG} explores 69\% (small) and 90\% (large) fewer labels than ASCG-ML and is 29\% and 31\% faster in mean time. The combination has the fewest iterations and the lowest mean time. Against ASCG-ML it needs 21\% and 13\% fewer iterations and 37\% and 43\% less time, and against Plain CG 26\% and 27\% fewer iterations and 32\% and 54\% less time. These results are consistent with the two methods being complementary.}

\begin{table}[!htbp]
\centering
\caption{\rev{Comparison and combination with ASCG-ML. Means and medians are over the instances on which all methods finished (289 small, 207 large).}}
\label{tab:ascg_results}
\resizebox{\textwidth}{!}{
\begin{tabular}{llrrrr}
\toprule
\textbf{Set} & \textbf{Method} & \textbf{Mean iterations} & \textbf{Mean labels} & \textbf{Mean time (s)} & \textbf{Median time (s)} \\
\midrule
Small & Plain CG & 37.3 & 125{,}277 & 2.73 & 0.077 \\
Small & ASCG-ML & 34.6 & 123{,}481 & 2.98 & 0.034 \\
Small & \textsc{Coupled-CG} & 33.8 & 37{,}854 & 2.13 & 0.07 \\
Small & \textsc{Coupled-CG} + ASCG-ML & 27.5 & 36{,}592 & 1.87 & 0.05 \\
\midrule
Large & Plain CG & 278.1 & 1{,}019{,}798 & 142.47 & 5.31 \\
Large & ASCG-ML & 231.1 & 933{,}473 & 115.68 & 4.71 \\
Large & \textsc{Coupled-CG} & 219.7 & 97{,}865 & 79.26 & 3.98 \\
Large & \textsc{Coupled-CG} + ASCG-ML & 201.8 & 98{,}301 & 65.98 & 3.01 \\
\bottomrule
\end{tabular}
}
\end{table}

\subsection{\rev{Comparison with Neural Routing Solvers}}
\label{appendix:neural_solver_comparison}

\rev{MTPOMO~\citep{liu2024multitask}, MVMoE~\citep{zhou2024mvmoe}, and RouteFinder~\citep{berto2024routefinder} are end-to-end constructive heuristics with no bound on solution quality. \textsc{Coupled-CG} accelerates an exact method: CG returns the LP lower bound of the set-partitioning formulation and, from the generated columns, an integer solution with a certified gap to that bound. We therefore compare on solution quality. Solve time is not comparable, because a neural solver only constructs a solution while CG also proves the bound, and CG-specific measures such as iterations and labels do not exist for neural solvers.}

\rev{All methods are run on our test split of 300 small and 300 large instances. Each instance is converted to the input convention of the neural solvers: coordinates, time windows, and service durations are scaled by a common factor of 100, and demands by the vehicle capacity. We use the released size-50 checkpoints with 8-fold augmentation, both zero-shot and after fine-tuning on our training instances for 50 epochs. Every method minimizes total travel distance, and costs are mapped back to the original scale. We report the gap to the CG LP lower bound, $(z-z_{\mathrm{LP}})/z_{\mathrm{LP}}$, on the instances where CG terminated exactly within the time limit (299 small, 185 large). This bound holds for every feasible solution, and every exact CG variant attains it, so the reference does not depend on which CG method computes it.}

\rev{Table~\ref{tab:neural_solver_comparison} reports the results. The best neural solver, fine-tuned RouteFinder, is on average 8.03\% (small) and 14.63\% (large) above the lower bound. The integer solutions of \textsc{Coupled-CG} are within 0.56\% and 1.87\%, and they reach the bound on 288 of 299 small and 133 of 185 large instances, against 248 and 88 for Plain CG. The mean gaps of the two CG variants differ by at most 0.10 percentage points. The released checkpoints have mean gaps of 16--36\% zero-shot and 8--23\% after fine-tuning on our training instances.}

\begin{table}[!htbp]
\centering
\caption{\rev{Gap to the CG LP lower bound on 299 small and 185 large instances. ``At bound'' counts the instances whose solution attains the bound. The CG rows use the integer solution obtained from the generated columns.}}
\label{tab:neural_solver_comparison}
\resizebox{\textwidth}{!}{
\begin{tabular}{lrrcrrc}
\toprule
 & \multicolumn{3}{c}{\textbf{Small}} & \multicolumn{3}{c}{\textbf{Large}} \\
\cmidrule(lr){2-4}\cmidrule(lr){5-7}
\textbf{Method} & \textbf{Mean gap} & \textbf{Median gap} & \textbf{At bound} & \textbf{Mean gap} & \textbf{Median gap} & \textbf{At bound} \\
\midrule
MTPOMO, zero-shot & 20.99\% & 15.62\% & 14/299 & 33.72\% & 17.27\% & 2/185 \\
MVMoE, zero-shot & 26.76\% & 21.10\% & 10/299 & 36.19\% & 22.76\% & 1/185 \\
RouteFinder, zero-shot & 16.39\% & 12.80\% & 12/299 & 23.52\% & 16.85\% & 0/185 \\
MTPOMO, fine-tuned & 11.29\% & 7.18\% & 43/299 & 18.83\% & 13.43\% & 2/185 \\
MVMoE, fine-tuned & 13.17\% & 9.65\% & 40/299 & 22.80\% & 18.55\% & 0/185 \\
RouteFinder, fine-tuned & 8.03\% & 5.83\% & 43/299 & 14.63\% & 11.55\% & 0/185 \\
\midrule
Plain CG & 0.66\% & 0.00\% & 248/299 & 1.86\% & 0.09\% & 88/185 \\
\textsc{Coupled-CG} & 0.56\% & 0.00\% & 288/299 & 1.87\% & 0.10\% & 133/185 \\
\bottomrule
\end{tabular}
}
\end{table}

%% file: Content/Conclusion.tex
\section{Conclusion}
\textsc{Coupled-CG} demonstrates that jointly coordinating pricing and column selection can yield additional gains beyond learning either decision independently. By propagating value across the two stages, the method reduces labeling effort and CG runtime while preserving exact termination, and its improvements extend to larger and structurally different CVRPTW instances.

Several limitations remain. First, the pricing policy is restricted to a predefined set of extension-ordering rules rather than making finer-grained label- or extension-level decisions, which limits the expressiveness of pricing control. Second, the learning objective uses cumulative labeling effort as a surrogate for computational cost. Although the experiments additionally evaluate end-to-end runtime, the training objective does not explicitly account for other computational components such as RMP reoptimization and neural inference. Third, \rev{although Section~\ref{appendix:generalization_beyond_cvrptw} tests two problems beyond CVRPTW, a systematic evaluation is still needed to identify which properties of a CG problem, such as its pricing structure and master-problem characteristics, determine when the coupled framework is effective.} Finally, the current study focuses on root-node CG; extending the framework to branch-and-price introduces additional interactions with branching decisions and changing pricing subproblems that are not considered here.

%% file: Content/Appendix.tex
\appendix

\section{Notation Table}\label{appendix:notation_tab}
\subsection{Notation}

Table~\ref{tab:notation} summarizes the main notation used throughout the
problem formulation and learning framework.
{\small
\begin{longtable}{p{0.19\linewidth} p{0.73\linewidth}}
\caption{Summary of the main notation.}\label{tab:notation}\\
\toprule
\textbf{Notation} & \textbf{Description} \\
\midrule
\endfirsthead
\toprule
\textbf{Notation} & \textbf{Description} \\
\midrule
\endhead

$\mathcal N$
& Set of customers. \\

$\mathcal R$
& Set of all feasible routes satisfying the routing constraints. \\

$\mathcal R_t$
& Set of routes currently included in the restricted master problem (RMP) at CG iteration $t$. \\

$c_r$
& Cost of route $r\in\mathcal R$. \\

$\alpha_{i,r}$
& Indicator equal to $1$ if route $r$ serves customer $i$, and $0$ otherwise. \\

$X_r$
& Master-problem decision variable associated with route $r$. \\

$V$
& Maximum number of available vehicles. \\

$d_i$
& Dual variable associated with the coverage constraint of customer $i$. \\

$q$
& Dual variable associated with the vehicle-limit constraint. \\

$\bar c_r$
& Reduced cost of route $r$. \\

\midrule

$I$
& Routing instance drawn from the target instance distribution $\mathcal D$. \\

$\mathcal D$
& Distribution of routing instances. \\

$t$
& Column-generation iteration index. \\

$T$
& Number of CG decision iterations before termination. \\

$P_t$
& Pricing problem at CG iteration $t$. \\

$N_t(g_t)$
& Number of labels explored during pricing at iteration $t$ under extension-ordering rule $g_t$. \\

$N_{\mathrm{scale}}$
& Fixed global normalization constant used to normalize labeling effort. \\

$Z_t$
& RMP objective value before the column insertion at iteration $t$. \\

$Z_1$
& Initial RMP objective value of the CG episode. \\

\midrule

$\mathcal G$
& Set of candidate extension-ordering rules. \\

\rev{$g_t$}
& \rev{Extension-ordering action selected at iteration $t$, with $g_t\in\mathcal G$.} \\

\rev{$\hat g_t$}
& \rev{Logged extension-ordering action observed in the offline dataset at iteration $t$.} \\

$\mathcal C_t$
& Candidate-column set returned by pricing at iteration $t$. \\

$M$
& Maximum number of candidate columns retained after pricing. \\

$k$
& Target number of columns selected for insertion at each CG iteration. \\

$O_t$
& Selected column set at iteration $t$, with
$O_t\subseteq\mathcal C_t$ and $|O_t|=k$. \\

$\pi^P$
& Extension-ordering policy. \\

$\pi^C$
& Column-selection policy. \\

\midrule

$s_t^P$
& Pre-pricing state defined by the current RMP at iteration $t$. \\

$s_t^C$
& Post-pricing state after applying $g_t$ and generating $\mathcal C_t$. \\

$h_t$
& Learned representation of the pre-pricing state,
$h_t=\operatorname{GNN}_{\theta}(s_t^P)$. \\

$\widetilde h_t$
& Learned post-pricing representation,
$\widetilde h_t=\operatorname{Enc}_{\theta}^{C}(h_t,\mathcal C_t)$. \\

$z_{t,i}$
& Relation-aware representation of candidate column $i\in\mathcal C_t$. \\

\midrule
$\xi_g$
& Learned prototype representation of extension-ordering rule $g$, used to estimate action support. \\

$\omega^P(h_t,g)$
& Prototype-based support score of extension-ordering action $g$ at state $h_t$. \\

$u^P(h_t,g)$
& Support uncertainty,
$u^P(h_t,g)=1-\omega^P(h_t,g)$. \\

$\kappa_P$
& Coefficient controlling the support-uncertainty penalty during extension-ordering selection. \\

$u_i$
& Unary value assigned to candidate column $i$. \\

$v_{ij}$
& Symmetric pairwise compatibility value between candidate columns $i$ and $j$. \\

$\Delta_{k}(i\mid O)$
& Marginal value of adding candidate $i$ to the current selected set $O$. \\

\midrule

$r_t^P$
& Pricing-stage reward based on normalized labeling effort. \\

$r_t^C$
& Column-selection reward based on the RMP objective improvement. \\

$\beta$
& Coefficient controlling the strength of the RMP-improvement shaping reward. \\

$Q_\psi^P(h_t,g)$
& Action-value function for extension-ordering action $g$. \\

$Q_{\phi,k}^C(\widetilde h_t,O)$
& Relation-conditioned value of selected column set $O$ with cardinality $k$. \\

$y_t^P$
& Temporal-difference target for the extension-ordering value function. \\

$y_t^C$
& Temporal-difference target for the column-selection value function. \\

\midrule

$\mathcal D_{\mathrm{off}}$
& Processed offline dataset of stage-specific CG transitions used for value learning. \\
$\mathrm{done}_t$
& Episode-termination indicator equal to $1$ if CG terminates after iteration $t$, and $0$ otherwise. \\

\bottomrule
\end{longtable}
}

\section{Proof of Proposition~\ref{prop:reward_consistency}}
\label{appendix:reward_consistency_proof}

\begin{proof}
For a CG trajectory with $T$ transitions, the cumulative pricing reward is
\begin{equation}
\sum_{t=1}^{T} r_t^P
=
-\frac{1}{N_{\mathrm{scale}}}
\sum_{t=1}^{T} N_t(g_t).
\end{equation}
Since $N_{\mathrm{scale}}>0$ is fixed across policies, maximizing this term is equivalent to minimizing cumulative labeling effort.

The cumulative column-selection reward telescopes as
\begin{equation}
\sum_{t=1}^{T} r_t^C
=
\frac{\beta}{\max\{1,|Z_1|\}}
\sum_{t=1}^{T}(Z_t-Z_{t+1})
=
\beta
\frac{Z_1-Z_{T+1}}
{\max\{1,|Z_1|\}}.
\end{equation}
All policies start from the same initial RMP, so $Z_1$ is fixed. Under exact CG termination, no negative-reduced-cost column remains; hence, $Z_{T+1}$ equals the optimal objective of the full LP relaxation and is also policy-independent. Therefore, the cumulative column-selection reward is constant with respect to the learned policies.

Combining the two reward components gives
\begin{equation}
\sum_{t=1}^{T}(r_t^P+r_t^C)
=
-\frac{1}{N_{\mathrm{scale}}}
\sum_{t=1}^{T}N_t(g_t)
+
\beta
\frac{Z_1-Z_{T+1}}
{\max\{1,|Z_1|\}}.
\end{equation}
The second term is policy-independent, yielding
\begin{equation}
\arg\max_{\pi^P,\pi^C}
\sum_{t=1}^{T}(r_t^P+r_t^C)
=
\arg\min_{\pi^P,\pi^C}
\sum_{t=1}^{T}N_t(g_t).
\end{equation}
Thus, under undiscounted return and exact termination, the RMP-improvement reward provides dense shaping without changing the optimizer of $\mathcal M$.
\end{proof}